\documentclass[a4paper,11pt]{article}
\usepackage[utf8]{inputenc}
\usepackage[english]{babel}
\usepackage{graphicx}
\usepackage{amsmath,amssymb, amsbsy, amsfonts, amsthm}
\usepackage[left=2.1cm,top=1.5cm,right=2.1cm, bottom=2.2cm,letterpaper]{geometry}
\usepackage{mathrsfs}
\usepackage{caption}
\usepackage{subcaption}
\usepackage{latexsym}
\usepackage{tcolorbox}
\usepackage{enumitem}
\usepackage{float}
\usepackage{hyperref}
\usepackage{url}
\usepackage[toc]{appendix}
\usepackage{xcolor}
\usepackage{color,psfrag}
\usepackage{etoc}
\usepackage{minitoc}
\usepackage{comment}
\usepackage{autonum}

\newtheorem{theorem}{Theorem}[section]
\newtheorem{corollary}[theorem]{Corollary}
\newtheorem{lemma}[theorem]{Lemma}
\newtheorem{definition}[theorem]{Definition}
\newtheorem{proposition}[theorem]{Proposition}

\newtheorem{remark}[theorem]{Remark}

\numberwithin{equation}{section}
\numberwithin{figure}{section}

\newcommand{\R}{\mathbb{R}}

\hypersetup{
	colorlinks=true,
	linkcolor=red,
	filecolor=magenta,
	urlcolor=cyan,
}

\makeatletter
\renewcommand*\env@matrix[1][*\c@MaxMatrixCols c]{%
	\hskip -\arraycolsep
	\let\@ifnextchar\new@ifnextchar
	\array{#1}}
\makeatother

\def\eq#1{(\ref{#1})}
\def\neweq#1{\begin{equation}\label{#1}}
\def\endeq{\end{equation}}

\begin{document}

\title{On the planar Taylor-Couette system and related exterior problems}

\author{Filippo GAZZOLA -- Ji\v{r}\'{\i} NEUSTUPA -- Gianmarco SPERONE}
\date{}
\maketitle
\vspace*{-6mm}

\begin{abstract}
	\noindent
We consider the planar Taylor-Couette system for the steady motion of a viscous incompressible fluid in the region between two concentric disks, the inner one being at rest and the outer one rotating with constant angular speed. We study the uniqueness and multiplicity of solutions to the forced system in different classes. For any angular velocity we prove that the classical Taylor-Couette flow is the unique smooth solution displaying rotational symmetry. Instead, we show that infinitely many solutions arise, even for arbitrarily small angular velocities, in the larger class of \textit{incomplete} solutions: roughly speaking, we constrain the test functions to belong
to a {\em smaller} space, that is, the solutions are allowed to belong to a {\em larger} space that includes possible singularities due to
the non-simply connectedness of the annulus. By prescribing the transversal flux, unique solvability of the Taylor-Couette system is recovered among rotationally-invariant incomplete solutions. Finally, we study the behavior of these solutions as the radius of the outer disk goes to infinity, connecting our results with the celebrated Stokes paradox.
\end{abstract}

{
	\hypersetup{linkcolor=black}
	\tableofcontents
}

\section{Introduction}

In his survey article on unsolved problems in Mathematical Analysis \cite[Problem 67]{maz2018seventy}, Maz'ya points out a classical but apparently forgotten issue about the stationary incompressible Navier-Stokes equations under Dirichlet boundary conditions in
{\em simply connected domains}: to show that uniqueness of the solution fails for large data.
In fact, multiplicity results were mostly obtained in the 1960's concerning \textit{open flows}, that is, the flow was assumed periodic in at least one direction: we recall the works of Velte \cite{velte1966stabilitat} and Yudovich \cite{judovich1966secondary} in 3D (the so-called \textit{Taylor-Couette problem}), and the article by Golovkin \cite{golovkin} in a 2D strip. One of the few available examples in the literature for a 3D bounded domain is also due to Yudovich \cite{yudovich1967example}. The significance of the works \cite{velte1966stabilitat,judovich1966secondary}, which deal with particular domains such as the (non-simply connected) region between two
concentric unbounded 3D cylinders, is that they settled a long-standing question left open since the famous experiments
of G.I.\ Taylor in 1923, see \cite{taylor1923}, illustrating the instability of the circular Couette flow. More recently, multiplicity
results for large Reynolds numbers in a (planar, simply connected) square for the stationary Navier-Stokes equations under
{\em Navier boundary conditions} have been obtained with computer assistance \cite{kochuniqueness}. Nevertheless, as far as our knowledge goes, the multiplicity question for the 2D Taylor-Couette problem has received little attention, and constitutes the main core of the present paper.\par
\begin{figure}[H]
	\begin{center}
	\includegraphics[scale=0.48]{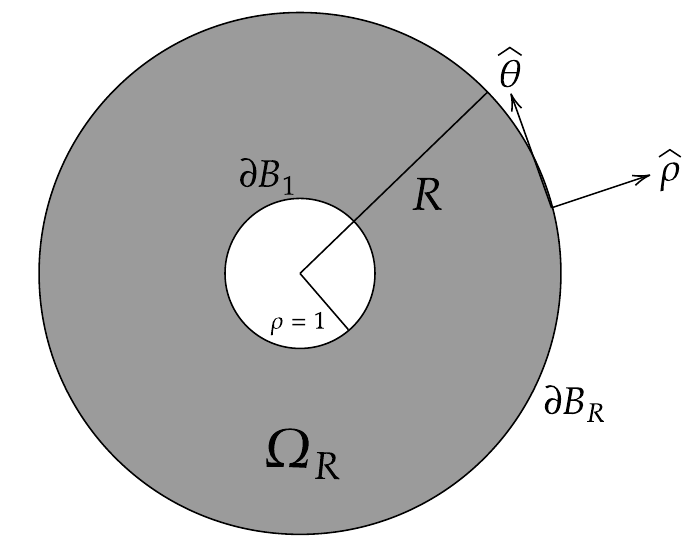}
	\end{center}
	\vspace*{-5mm}
	\caption{The annulus $\Omega_{R}$ in \eqref{annulus}.}\label{dom1}
\end{figure}

For any $R>1$, define the annulus $\Omega_R \subset \mathbb{R}^{2}$ as
\neweq{annulus}
\Omega_{R} = B_{R} \setminus \overline{B_{1}} \, ,
\endeq
with $B_{r}$ being the open disk of radius $r>0$ centered at the origin, see Figure \ref{dom1}.
Within the polar coordinate system $(\rho,\theta) \in [0,\infty) \times [0,2\pi)$, let $\lbrace \widehat{\rho}, \widehat{\theta} \rbrace \subset \mathbb{R}^{2}$ be the usual orthonormal basis, namely
\begin{equation}\label{rohat}
\widehat{\rho}=(\cos(\theta),\sin(\theta)) \quad \text{and} \quad \widehat{\theta}=(-\sin(\theta),\cos(\theta)) \qquad \forall \theta \in [0,
2\pi) \, ,
\end{equation}
see again Figure \ref{dom1}, and any point $\xi\in\R^2$ is denoted by $\xi = \rho\widehat{\rho}$.\par
We consider the following boundary-value problem associated with the steady-state Navier-Stokes equations in $\Omega_{R}$
(the kinematic viscosity has been set equal to 1):
\begin{equation}\label{nsstokes0}
\left\{
\begin{aligned}
& -\Delta u+(u\cdot\nabla)u+\nabla p=f \, ,\ \quad  \nabla\cdot u=0 \ \ \mbox{ in } \ \ \Omega_{R} \, , \\[3pt]
& u = \omega \, \widehat{\theta} \ \text{ on } \ \partial B_R \, , \qquad u = (0,0) \ \text{ on } \ \partial B_{1} \, .
\end{aligned}
\right.
\end{equation}
In \eqref{nsstokes0}, $u : \Omega_{R} \longrightarrow \mathbb{R}^2$ is the velocity vector field, $p : \Omega_{R} \longrightarrow \mathbb{R}$ is
the scalar pressure and $f : \Omega_{R} \longrightarrow \mathbb{R}^2$ is an external force acting on the fluid. The boundary conditions in
\eqref{nsstokes0}$_2$ dictate that the inner disk $\partial B_1$ remains at rest, while the outer disk $\partial B_R$ rotates with constant
angular speed $\omega \geq 0$. Notice that the compatibility condition
\begin{equation} \label{zeroflux2}
\int_{\Omega_R}\nabla \cdot u =\int_{\partial\Omega_R}u \cdot \nu = \omega\int_{\partial B_R} \widehat{\theta} \cdot \widehat{\rho} = 0
\end{equation}
is satisfied for the solutions to \eq{nsstokes0}. In  \eqref{zeroflux2}, $\nu$ denotes the outward unit normal to $\Omega_R$, so that
$\nu=-\widehat{\rho}$
on $\partial B_1$ and $\nu=\widehat{\rho}$ on $\partial B_R$. In the unforced case (when $f = 0$), an explicit classical solution to
\eq{nsstokes0} exists for every $\omega \geq 0$, and is given by
\begin{equation} \label{tc1}
v_{0}(\xi) := \dfrac{R \omega}{R^2 - 1} \left( \rho - \dfrac{1}{\rho} \right) \widehat{\theta} \quad \text{and} \quad q_{0}(\xi) :=
\dfrac{1}{2} \left( \dfrac{R \omega}{R^2 - 1} \right)^{2} \left( \rho^{2} - \dfrac{1}{\rho^{2}} - 4 \log(\rho) \right) \qquad \forall \xi \in
\Omega_{R} \, ,
\end{equation}
known in the literature as the \textit{Taylor-Couette flow} \cite[Chapter II]{landau}.\par
System \eq{nsstokes0} is the forced 2D version of the Taylor-Couette problem in a non-simply connected bounded
domain and uniqueness/multiplicity are related to the magnitude of the force $f$ and of the angular velocity $\omega$. Since it may be
$f\neq0$, we call the related problem to \eq{nsstokes0} the {\it generalized Taylor-Couette problem}, and its solution the {\it generalized Taylor-Couette flow}, see \eq{gtc1}.\par
We discuss the uniqueness and multiplicity of solutions to \eq{nsstokes0} from several points of view. In contrast with the 3D case \cite{velte1966stabilitat,judovich1966secondary}, in Section \ref{uniquetc} we prove with elementary methods that, for arbitrarily large data, the generalized
Taylor-Couette flow is the unique strong solution to \eq{nsstokes0} displaying rotational symmetry, see Theorem \ref{UVconstant0} and Corollary \ref{UVconstant0cor}. Instead, in Section \ref{singsolutions} we show that infinitely many solutions do exist, even for arbitrarily small data,
in a new and larger class of \textit{incomplete} solutions, see Definition \ref{incompletedef} and Theorem \ref{nonuniquetheomain}. This class is larger because the space of admissible test functions is taken smaller and, hence, they {\em are not weak solutions}
in the usual sense. Their mathematical interest is the connection between the nontrivial topology of $\Omega_R$ and the possible appearance of
singularities allowed by the larger space. In turn, their physical interest is that these singularities are due to a {\em pressure drop}, that is,
a discontinuity of the pressure that, apparently, does not give discontinuity of its gradient, see \eqref{firstchar} below. Exploiting the non-simply connectedness of $\Omega_R$ and following the ideas of Foias \& Temam \cite{foiastemam78}, the admissible test functions for the class of incomplete solutions have zero flux across each transversal section of $\Omega_{R}$. Unique solvability of \eq{nsstokes0} is then ensured
also among rotationally invariant incomplete solutions provided that their transversal flux is prescribed, see Theorem \ref{fluxtheo1}.\par
On the other hand, without the rotational symmetry assumption, unique solvability of \eq{nsstokes0}, in both classes of strong
and incomplete solutions, is guaranteed only under smallness conditions on the data of the problem (angular speed, external force,
transversal flux),
which we quantify. Several Poincaré-Sobolev constants are used in literature to ensure uniqueness to inhomogeneous Dirichlet problems for forced Navier-Stokes
equations, see \cite[Section IX.4]{galdi2011introduction}. Most of these statements merely give rough qualitative bounds on these constants.
Recently, in \cite{gazspefra,gazspe,cortona,gazzola2022connection} some {\em quantitative} (although not optimal) estimates for these constants in
non-simply connected domains were found. This technical part is postponed to Appendix \ref{bondsSob}, where we give
quantitative bounds for these constants in the annulus $\Omega_R$.\par
Finally, by studying the behavior of the solutions to \eqref{nsstokes0} as $R\to\infty$, we highlight a connection with the celebrated {\em Stokes paradox} \cite{stokes1851effect}. The classical result by Chang \& Finn \cite{chang1961solutions} shows that the linear version of \eqref{nsstokes0} with prescribed angular velocity at infinity has no solutions. In line with
other problems of hydrodynamics in planar exterior domains \cite{korobkov2019convergence,korobkov2020steady,korobkov2023stationary}, we raise the question on the solvability of the full nonlinear problem \eqref{nsstokes0} in $\mathbb{R}^2 \setminus \overline{B_{1}}$, see \eqref{nsexterior}:
in Theorem \eqref{nosimetria} we show that it admits no rotationally invariant solutions, regardless of the angular velocity at infinity.
In Remark \ref{NSparadox} we leave unanswered the possibility of a \textit{Navier-Stokes paradox}. In the same spirit, in Section \ref{stokesparadox} we study two related exterior Stokes problems and we propose some open questions on the existence of solutions with given growth rate at infinity.

\section{Solutions for the generalized Taylor-Couette problem} \label{uniquetc}

For our purposes, we need to consider both the Sobolev space $H_{0}^{1}(\Omega_{R})$ and the larger space of scalar functions vanishing only on
$\partial B_{1}$:
$$
H^1_*(\Omega_{R}) := \{v\in H^1(\Omega_{R}) \ | \ v=0 \ \ \mbox{on} \ \ \partial B_{1}\}\, .
$$
This space is the closure of $\mathcal{C}^{\infty}_{0}(\overline{B_{R}}\setminus\overline{B_{1}})$ with respect to the norm $v\mapsto\|\nabla
v\|_{L^2(\Omega_{R})}$: since $|\partial B_{1}|_{1}>0$ (the 1D-Hausdorff measure),
the Poincaré inequality holds in $H^1_*(\Omega_{R})$, showing that $v\mapsto\|\nabla v\|_{L^2(\Omega_{R})}$ is indeed a norm on $H^1_*(\Omega_{R})$. We also consider the two functional spaces
$$
H^{1}_{*,\sigma}(\Omega_{R}) = \{ v \in H^1_*(\Omega_{R}) \ | \ \nabla \cdot v = 0 \ \ \text{in} \ \ \Omega_{R} \} \quad \text{and} \quad
H^{1}_{0,\sigma}(\Omega_{R}) = \{ v \in H^{1}_{0}(\Omega_{R}) \ | \ \nabla \cdot v = 0 \ \ \text{in} \ \ \Omega_{R} \} \, .
$$
They are Hilbert spaces if endowed with the $L^{2}(\Omega_{R})$-scalar product of the gradients.\par
Given $f \in L^{2}(\Omega_{R})$, a vector field $u\in H^{1}_{*,\sigma}(\Omega_{R})$ is called a \textit{weak solution} to \eqref{nsstokes0} if $u$
verifies \eqref{nsstokes0}$_2$ in the trace sense and
\begin{equation} \label{nstokesdebil}
\int_{\Omega_{R}} \nabla u \cdot \nabla \varphi + \int_{\Omega_{R}} (u \cdot \nabla) u \cdot \varphi  = \int_{\Omega_{R}} f \cdot \varphi \qquad
\forall \varphi \in H^{1}_{0,\sigma}(\Omega_{R}) \, .
\end{equation}	
Since $\partial \Omega_{R}$ and the boundary data in \eqref{nsstokes0}$_2$ are of class $\mathcal{C}^{\infty}$, well-known regularity results
\cite[Theorem IX.5.2]{galdi2011introduction} imply that any weak solution satisfies $u \in H^{2}(\Omega_{R})$ and that there exists an
associated pressure
$p \in H^{1}(\Omega_{R})$ such that the pair $(u,p)$ solves \eqref{nsstokes0} in strong form. Therefore, since we always assume that $f\in
L^2(\Omega_R)$,
in the sequel we make no distinction between weak and strong solutions to \eqref{nsstokes0}.\par
Next, we observe that the functions in \eqref{tc1} are \textit{rotationally invariant}, meaning that
\begin{equation} \label{rotinvariant}
v_{0}(\xi) = \mathcal{R}(\phi)^{\top} v_{0}(\mathcal{R}(\phi) \xi) \qquad \text{and} \qquad q_{0}(\xi) = q_{0}(\mathcal{R}(\phi) \xi) \qquad
\forall \xi \in \Omega_{R} \, , \ \forall \phi \in [0, 2\pi] \, ,
\end{equation}
where
\begin{equation} \label{rotmat}
\mathcal{R}(\phi) :=
\begin{bmatrix}
	\cos(\phi)       & -\sin(\phi)  \\[0.2cm]
	\sin(\phi)       & \cos(\phi)
\end{bmatrix}
\end{equation}
is the rotation matrix about the origin by an angle $\phi$. A rotationally invariant vector field $f \in L^2(\Omega_{R})$ can be expressed as
\begin{equation}\label{frotational}
f(\xi) = f^{\rho}(\rho) \widehat{\rho} + f^{\theta}(\rho) \widehat{\theta} \qquad \mbox{for a.e. }\xi\in\Omega_{R}\, ,
\end{equation}
for some scalar functions $f^{\rho}, f^{\theta}  \in L^2(1,R)$. Then
\begin{equation}\label{L2}
\|f\|_{L^2(\Omega_R)}^2=\|f^\rho\|_{L^2(\Omega_R)}^2+\|f^\theta\|_{L^2(\Omega_R)}^2=2\pi\int_1^{R} \rho
\left(f^\rho(\rho)^2+f^\theta(\rho)^2\right)\, d\rho\, .
\end{equation}
Henceforth, given any subspace $X \subseteq L^{2}(\Omega_{R})$ (of vector fields or scalar functions), we will denote by
\begin{center}
$\mathcal{R}[X]$ the
subspace of $X$ comprising rotationally invariant functions.
\end{center}
Let
$$
A_{R}(\omega,f) := \dfrac{R}{R^2 - 1} \left(\omega + \dfrac{R}{2}\int_{1}^{R} f^{\theta}(t)\, dt - \dfrac{1}{2R}\int_{1}^{R}
t^{2}f^{\theta}(t)\, dt \right) \,,
$$
and let $\lambda(\Omega_{R})>0$ be the least Laplace-Dirichlet eigenvalue in $\Omega_R$, see \eqref{eigen1}. We then prove

\begin{theorem}\label{UVconstant0}
Given $\omega\ge0$ and $f \in \mathcal{R}[L^{2}(\Omega_{R})]$, the generalized Taylor-Couette flow defined in $\Omega_{R}$ by
\begin{equation} \label{gtc1}
\left\{
\begin{aligned}
& v_{*}(\xi) := \left[ A_{R}(\omega,f) \left( \rho - \dfrac{1}{\rho} \right) - \dfrac{\rho}{2} \int_{1}^{\rho} f^{\theta}(t) \, dt +
\dfrac{1}{2\rho} \int_{1}^{\rho} t^{2} f^{\theta}(t) \, dt \right] \widehat{\theta} \\[6pt]
& q_{*}(\xi) := \int_{1}^{\rho} \dfrac{1}{t} \left[ A_{R}(\omega,f) \left( t - \dfrac{1}{t} \right) - \dfrac{t}{2} \int_{1}^{t} f^{\theta}(s)
\, ds + \dfrac{1}{2t} \int_{1}^{t} s^{2} f^{\theta}(s) \, ds \right]^{2} dt +  \int_{1}^{\rho} f^{\rho}(t) \, dt
\end{aligned}
\right.
\end{equation}
solves \eqref{nsstokes0}. Moreover, \eqref{gtc1} is the unique strong (and weak) solution to \eqref{nsstokes0}
\begin{itemize}
\item in the space $\mathcal{R}[H^{2}(\Omega_{R}) \cap H^{1}_{*,\sigma}(\Omega_{R}) ] \times \mathcal{R}[H^{1}(\Omega_{R})]/\mathbb{R}$ for any $\omega \geq 0$ and $f \in \mathcal{R}[L^{2}(\Omega_{R})]$;

\item in the whole space $[H^{2}(\Omega_{R}) \cap H^{1}_{*,\sigma}(\Omega_{R})] \times H^{1}(\Omega_{R})/\mathbb{R}$ whenever
\begin{equation} \label{umbral1cor}
K_R(\omega,f) := \omega+\frac{\sqrt{R-1}}{\sqrt{15}R}\sqrt{8R^4-7R^3-7R^2+3R+3}\, \|f^\theta\|_{L^2(\Omega_R)}<\sqrt{\lambda(\Omega_R)}\, .
\end{equation}
\end{itemize}
\end{theorem}
\noindent
\begin{proof} Let $(v,p) \in \mathcal{R}[H^{2}(\Omega_{R})] \times \mathcal{R}[H^{1}(\Omega_{R})]$ be a strong solution to \eqref{nsstokes0}.
Then, $p$ depends only on $\rho \in [1,R]$, while the velocity field and external force may be written as
$$
v(\xi) = v^{\rho}(\rho) \widehat{\rho} + v^{\theta}(\rho) \widehat{\theta} \quad \text{and} \quad f(\xi) = f^{\rho}(\rho) \widehat{\rho} +
f^{\theta}(\rho) \widehat{\theta} \qquad \text{for a.e. }\xi \in \Omega_{R} \, ,
$$
for some functions $v^{\rho}, v^{\theta}  \in H^{2}(1,R)$ and $f^{\rho}, f^{\theta}  \in L^{2}(1,R)$. Since $v^\theta$
is independent
of $\theta$, the incompressibility condition becomes
$$
\dfrac{d}{d\rho}(\rho v^{\rho}(\rho))=0\quad \Longrightarrow \quad v^{\rho}(\rho) = \dfrac{A}{\rho} \qquad\text{for a.e. }\rho\in(1,R) \, ,
$$
for some constant $A \in \mathbb{R}$. The boundary conditions \eqref{nsstokes0}$_2$ imply $A=0$, so that $v(\xi) = v^{\theta}(\rho)
\widehat{\theta}$ for a.e.\
$\xi \in \Omega_{R}$, where $v^{\theta}(1) = 0$ and $v^{\theta}(R) = \omega$. Then, the first equation in \eqref{nsstokes0}$_1$ reduces to
$$
\frac{d p}{d \rho}(\rho) = \dfrac{1}{\rho} v^{\theta}(\rho)^{2} + f^{\rho}(\rho) \quad \text{and} \quad -\dfrac{d^2 v^{\theta}}{d\rho^2}(\rho) -
\frac{1}{\rho} \dfrac{d v^{\theta}}{d\rho}(\rho) + \dfrac{1}{\rho^{2}} v^{\theta}(\rho) = f^{\theta}(\rho) \qquad \text{for a.e.} \ \rho \in (1,R)
\, ,
$$
thereby yielding $(v,p)\equiv(v_{*},q_{*})$, as defined in \eqref{gtc1}. Not only this proves that \eqref{gtc1} solves \eqref{nsstokes0}
(which could have been obtained by direct computations), but also that it is the unique solution within the class of rotationally
invariant solutions. We have so proved the first item, namely the unique solvability of \eqref{nsstokes0} in the class of rotationally
invariant solutions independently of the size of the data $\omega$ and $f$.
\par
Moreover, the explicit form \eqref{gtc1} enables us to bound the maximum modulus of $v_{*}$. Indeed,
$$
\begin{aligned}
|v_{*}(\xi)| &= \left|A_{R}(\omega,f) \left(\rho-\dfrac{1}{\rho}\right)-\dfrac{\rho}{2} \int_{1}^{\rho} f^{\theta}(t) \, dt + \dfrac{1}{2\rho}
\int_{1}^{\rho} t^{2} f^{\theta}(t) \, dt\right|\\[6pt]
 &\le \dfrac{R}{R^2-1}\dfrac{\rho^2-1}{\rho}\left[\omega+\frac{1}{2R}\int_{1}^{R}(R^2-t^2)|f^{\theta}(t)|\, dt\right]
+\frac1{2\rho} \int_{1}^{\rho}(\rho^2-t^2)|f^{\theta}(t)|\, dt\\[6pt]
   &\le \omega+\frac{1}{2R}\int_{1}^{R}(R^2-t^2)|f^{\theta}(t)|\, dt+\frac1{2\rho} \int_{1}^{\rho}(\rho^2-t^2)|f^{\theta}(t)|\, dt\\[6pt]
   &\le \omega+\frac{1}{2R}\left[\int_{1}^{R}(R^2-t^2)^2\, dt\right]^{1/2}\|f^{\theta}\|_{L^2(1,R)}
   +\frac1{2\rho}\left[\int_{1}^{\rho}(\rho^2-t^2)^2\, dt\right]\|f^{\theta}\|_{L^2(1,\rho)}\\[6pt]
   &\le \omega+\left[\frac{\sqrt{R-1}}{2\sqrt{15}R}\sqrt{8R^4-7R^3-7R^2+3R+3}
   +\frac{\sqrt{\rho-1}}{2\sqrt{15}\rho}\sqrt{8\rho^4-7\rho^3-7\rho^2+3\rho+3}\right]\|f^{\theta}\|_{L^2(1,R)}\\[6pt]
   &\le \omega+\frac{\sqrt{R-1}}{\sqrt{15}R}\sqrt{8R^4-7R^3-7R^2+3R+3}\, \|f^{\theta}\|_{L^2(\Omega_R)} \, ,
\end{aligned}
$$
which implies (see \eqref{umbral1cor}) that
\begin{equation}\label{Linfty}
\|v_{*}\|_{L^{\infty}(\Omega_{R})}\le K_{R}(\omega,f)\, .
\end{equation}

The second item is a variant of a well-known result, see \cite[Section IX.4]{galdi2011introduction}: unique solvability for \eqref{nsstokes0}
under a
smallness condition on the boundary velocity $\omega$ and on the force $f$. The novelty here is the explicit bound involving the first eigenvalue
of
the Laplace-Dirichlet operator in $\Omega_{R}$, \eqref{eigen1}. In order to prove it, suppose that $u \in H^{1}_{*,\sigma}(\Omega_{R})$ is another strong
solution to  \eqref{nsstokes0}. Let $z := u - v_{*} \in H^{1}_{0,\sigma}(\Omega_{R})$ and subtract the equations \eqref{nstokesdebil}
corresponding to $u$ and $v_{*}$, thereby obtaining
$$
\int_{\Omega_{R}} \nabla z \cdot \nabla \varphi + \int_{\Omega_{R}} (u \cdot \nabla) z \cdot \varphi + \int_{\Omega_{R}} (z \cdot \nabla) v_{*}
\cdot \varphi = 0 \qquad \forall \varphi \in H^{1}_{0,\sigma}(\Omega_{R}) \, .
$$
By taking $\varphi=z$ and noticing (after an integration by parts) that
$$
\int_{\Omega_{R}} (u \cdot \nabla) z \cdot z = 0 \qquad \text{and} \qquad \int_{\Omega_{R}} (z \cdot \nabla) v_{*} \cdot z = - \int_{\Omega_{R}}
(z \cdot \nabla) z \cdot v_{*} \, ,
$$
we deduce, from the H\"older inequality, from \eqref{sobpoin}, and from \eqref{Linfty}, that
$$
\|\nabla z\|^2_{L^2(\Omega_{R})} = \int_{\Omega_{R}} (z \cdot \nabla) z \cdot v_{*}\leq \|\nabla z\|_{L^2(\Omega_{R})} \|z\|_{L^2(\Omega_{R})} \|
v_{*} \|_{L^{\infty}(\Omega_{R})}\leq \dfrac{K_{R}(\omega,f)}{\sqrt{\lambda(\Omega_R)}}\|\nabla z\|^{2}_{L^2(\Omega_{R})}\, ,
$$
so that unique solvability for \eqref{nsstokes0} is ensured whenever \eqref{umbral1cor} holds.\end{proof}

A first straightforward consequence of Theorem \ref{UVconstant0} is

\begin{corollary}\label{UVconstant0cor}
If $f=0$, the Taylor-Couette flow \eqref{tc1} is the unique strong (classical) solution to \eqref{nsstokes0}
\begin{itemize}
\item in the space $\mathcal{R}[H^{2}(\Omega_{R}) \cap H^{1}_{*,\sigma}(\Omega_{R})] \times \mathcal{R}[H^{1}(\Omega_{R})]/\mathbb{R}$ for any $\omega \geq 0$;
\item in the whole space $[H^{2}(\Omega_{R}) \cap H^{1}_{*,\sigma}(\Omega_{R})]\times H^{1}(\Omega_{R})/\mathbb{R}$ whenever $\omega < \sqrt{\lambda(\Omega_{R})}$.
\end{itemize}
\end{corollary}

A second consequence of Theorem \ref{UVconstant0} shows that the fluid remains at rest under the action of an arbitrarily
large radial and rotationally invariant force.

\begin{corollary} \label{atrest}
Let $\omega \geq 0$ and $f \in \mathcal{R}[L^{2}(\Omega_{R})]$.
\begin{itemize}	
\item If $f^{\theta}=0$ the velocity field $v_{*}$ in \eqref{gtc1} coincides with $v_{0}$ in \eqref{tc1}.
\item If $f^{\theta}=0$ and $\omega=0$, the generalized Taylor-Couette flow \eqref{gtc1} becomes
$$
v_{*}(\xi) = 0 \qquad \text{and} \qquad q_{*}(\xi) =  \int_{1}^{\rho} f^{\rho}(t) \, dt \qquad \forall \xi \in \overline{\Omega_{R}} \, .
$$
\end{itemize}
\end{corollary}

Theorem \ref{UVconstant0} and Corollaries \ref{UVconstant0cor}-\ref{atrest} deserve some comments.

\begin{remark}
The inequality \eqref{umbral1cor} or rotational invariance are sufficient conditions for the unique weak solvability
of \eqref{nsstokes0} which is guaranteed even for an arbitrarily large radial force, since the inequality in \eqref{umbral1cor} only involves
the angular component of the force; see Corollary \ref{atrest}. Related Liouville-type results for a generalized 3D Taylor-Couette problem have been recently published by Kozono et al. in \cite{kozono2023liouville}.
When $f \in \mathcal{R}[L^{2}(\Omega_{R})]$, in the next section we show that multiplicity arises in a wider class of solutions.
\end{remark}

\begin{remark} \label{infinitas}
The violation of \eqref{umbral1cor} may only lead to the existence of a non-rotationally invariant strong solution to \eqref{nsstokes0}.
Furthermore, notice that if
$(u,p) \in H^2(\Omega_{R})\times H^1(\Omega_{R})$ is a strong solution to \eqref{nsstokes0} with a given rotationally invariant external force
\eqref{frotational}, then the pair
\begin{equation}\label{rotinv}
u^{\phi}(\xi) := \mathcal{R}(\phi)^{\top} u(\mathcal{R}(\phi) \xi) \qquad \text{and} \qquad p^{\phi}(\xi) := p(\mathcal{R}(\phi) \xi)
\qquad \forall \xi \in \Omega_{R} \, ,
\end{equation}
is also a solution to \eqref{nsstokes0}, for every $\phi \in [0,2 \pi)$. This means that, if a bifurcation occurs, then infinitely many
non-rotationally invariant strong solutions to \eqref{nsstokes0} exist.
\end{remark}

\begin{remark}
In favor of a possible unconditional uniqueness in the unforced case, notice that
	\begin{equation}\label{nsstokes0br}
	\left\{
	\begin{aligned}
	& -\Delta u+(u\cdot\nabla)u+\nabla p=0 \, ,\ \quad  \nabla\cdot u=0 \ \ \mbox{ in } \ \ B_{R} \, , \\[3pt]
	& u = \omega \, \widehat{\theta} \ \text{ on } \ \partial B_R \, ,
	\end{aligned}
	\right.
	\end{equation}
	has a unique classical solution for every $\omega \geq 0$ and $R>0$, given explicitly by
	$$
	\bar{u}(\xi) := \dfrac{\omega \rho}{R} \, \widehat{\theta} \qquad \text{and} \qquad \bar{p}(\xi) := \dfrac{1}{2} \left( \dfrac{\omega \rho}{R} \right)^{2} \qquad \forall \xi \in B_{R} \, .
	$$
This follows from the fact that $\rho \, \widehat{\theta}=(-y,x)$ and the identity $(z \cdot \nabla) \bar{u} \cdot z \equiv 0$ in $B_{R}$, for every $z \in \mathbb{R}^{2}$.
\end{remark}

Corollary \ref{UVconstant0} ensures unconditional uniqueness
of weak solutions to \eqref{nsstokes0} only if $\omega>0$ is sufficiently small (for a fixed $R>1$), or if $R > 1$ is sufficiently small
(for a fixed $\omega>0$), see \eqref{lambda}. This is why we conclude this section by investigating what happens as $R\to\infty$, including the case $R=\infty$.
This is related to the celebrated \textit{invading domains} technique by Jean Leray \cite{leray1933etude} that we will analyze for the Stokes equations in Section \ref{stokesparadox}.
\par
Assuming that $f=0$ and fixing $\omega>0$, Corollary \ref{UVconstant0cor} states that the unique weak (classical) rotationally
invariant solution to \eqref{nsstokes0} is given by the Taylor-Couette flow \eqref{tc1}, namely
$$
v_R(\xi) := \dfrac{R \omega}{R^2 - 1} \left( \rho - \dfrac{1}{\rho} \right) \widehat{\theta} \quad \text{and} \quad q_R(\xi)=
\dfrac{1}{2}\left(\dfrac{R \omega}{R^2 - 1}\right)^{2}\left(\rho^{2}-\dfrac{1}{\rho^{2}}-4\log(\rho)\right)\qquad\forall\xi\in\Omega_{R}\, ,
$$
where we emphasized the dependence on $R$. It is straightforward that, if we extend $v_R$ and $q_R$ by
$$v_R(\xi) = \omega\widehat{\theta}\, ,\qquad q_R(\xi) = \dfrac{1}{2} \left( \dfrac{R \omega}{R^2 - 1} \right)^{2}
\left(R^{2}-\dfrac{1}{R^{2}}-4\log(R)\right)\qquad \forall \xi \in \R^2 \setminus B_{R} \, ,$$
we obtain continuous functions defined over the whole exterior domain $\Omega_{\infty} := \R^2\setminus B_1$ which satisfy
$$
v_R, \, q_R\to0\quad\mbox{in } \ L^\infty_{\rm loc}(\Omega_\infty)\quad\mbox{as } \ R\to\infty\, .
$$

On the other hand, one may take directly $R=\infty$, without approximating through an invading domains procedure. We discuss
whether there exists a classical solution $(v,q)\in\mathcal{C}^{2}(\overline{\Omega_{\infty}}) \times \mathcal{C}^{1}(\overline{\Omega_{\infty}})$
to the following exterior problem:
\begin{equation}\label{nsexterior}
\left\{
\begin{aligned}
	& -\Delta v+ (v \cdot \nabla)v +  \nabla q=0 \, ,\ \quad  \nabla\cdot v=0 \ \ \mbox{ in } \ \ \Omega_{\infty}\, ,
\\[3pt]
	& \lim\limits_{\rho \to \infty} v(\xi) = \omega \, \widehat{\theta}\, , \qquad v = (0,0) \ \text{ on } \ \partial B_{1} \, .
\end{aligned}
\right.
\end{equation}
Notice that \eqref{nsexterior} admits no solutions having a finite Dirichlet integral, since
$$
| \nabla \widehat{\theta}(\xi)|=\dfrac{1}{\rho}\qquad\forall \xi \in \Omega_{\infty} \, .
$$
Hence, for every solution $(v,q)$ to \eqref{nsexterior}, we have that $\nabla v \notin L^{2}(\Omega_{\infty})$. Moreover, the
following holds:

\begin{theorem} \label{nosimetria}
Given $\omega > 0$, the exterior problem \eqref{nsexterior} admits no rotationally invariant solutions.
\end{theorem}
\noindent
\begin{proof} Let $(v, q) \in \mathcal{C}^{2}(\overline{\Omega_{\infty}}) \times \mathcal{C}^{1}(\overline{\Omega_{\infty}})$ be a rotationally
invariant solution to \eqref{nsexterior}. Then, the pressure $q$ depends only on $\rho \in [1,+\infty)$, while the velocity
field may be written as
$$
v(\xi) = v^{\rho}(\rho) \widehat{\rho} + v^{\theta}(\rho) \widehat{\theta} \qquad  \forall \xi \in \Omega_{\infty} \, ,
$$
for some $v^{\rho}, v^{\theta}  \in \mathcal{C}^{\infty}(1,+\infty)$. Since $v^\theta$ is independent
of $\theta$, the incompressibility condition becomes
$$
\dfrac{d}{d\rho}\big(\rho v^{\rho}(\rho)\big)=0\quad \Longrightarrow \quad v^{\rho}(\rho) = \dfrac{A}{\rho} \qquad \forall \rho > 1 \, ,
$$
for some constant $A \in \mathbb{R}$. The boundary conditions \eqref{nsexterior}$_2$ imply $A=0$, so that $v(\xi) = v^{\theta}(\rho)
\widehat{\theta}$ for every
$\xi \in \Omega_{\infty}$, with $v^{\theta}(1) = 0$ and $v^{\theta}\to \omega$ as $\rho \to +\infty$. Then, the first equation in
\eqref{nsexterior}$_1$ reduces to
$$
\frac{d p}{d \rho}(\rho) = \dfrac{1}{\rho} v^{\theta}(\rho)^{2} \quad \text{and} \quad -\dfrac{d^2 v^{\theta}}{d\rho^2}(\rho) - \frac{1}{\rho}
\dfrac{d v^{\theta}}{d\rho}(\rho) + \dfrac{1}{\rho^{2}} v^{\theta}(\rho) = 0 \qquad \forall \ \rho > 1 \, ,
$$
thereby yielding
$$
v^{\theta}(\rho) = \dfrac{B}{\rho} + C \rho \qquad \forall \rho \geq 1 \, , \quad \text{for some constants} \  B,C \in \mathbb{R} \, .
$$
Enforcing the boundary conditions \eqref{nsexterior}$_2$ gives $B=C=0$, a contradiction.
\end{proof}

\begin{remark}\label{NSparadox}
In view of Remark \ref{infinitas}, as a consequence of Theorem \ref{nosimetria} we deduce that \eqref{nsexterior} is either not
solvable, or it has infinitely many non-rotationally invariant solutions. Furthermore, with the same proof
as Theorem \ref{nosimetria}, we deduce that the exterior problem \eqref{nsexterior} {\em without condition at infinity} only admits
as rotationally invariant solutions the family of functions
\neweq{uniqueK}
v(\xi)=K\left(\rho-\frac{1}{\rho}\right)\widehat\theta\, ,\quad q(\xi)= \dfrac{K^2}{2} \left( \rho^2 -\frac{1}{\rho^2} - 4\log(\rho) \right)
\quad\forall\xi\in\Omega_\infty\qquad(K\in\mathbb{R})\, .
\endeq
\end{remark}

\section{Non-uniqueness in a larger class of solutions} \label{singsolutions}

Theorem \ref{UVconstant0} does not exclude multiplicity of weak solutions to \eqref{nsstokes0}
(in fact, strong solutions since we assume that $f\in L^2(\Omega_R)$) in the class of non-rotationally invariant
functions if \eqref{umbral1cor} is violated. In view of \eqref{lambda}, the latter occurs if $\omega>0$ is
sufficiently large (for a fixed $R>1$), or if $R > 1$ is sufficiently small (for a fixed $\omega>0$). But Theorem \ref{UVconstant0} also
suggests to investigate whether multiplicity results for \eqref{nsstokes0} in $H^{1}_{*,\sigma}(\Omega_{R})$ can be proved with a
different notion of solution, e.g., by restricting the class $H^{1}_{0,\sigma}(\Omega_{R})$ of test functions in \eqref{nstokesdebil}.
To this end, we need a refinement of the classical Helmholtz-Weyl decomposition for vector fields \cite{Helmholtz}, which states that
\begin{equation}\label{HW}
L^2(\Omega_R)=G_{1}(\Omega_R) \oplus G_{2}(\Omega_R) \, ,
\end{equation}
where
$$
\begin{aligned}
& G_{1}(\Omega_R) := \{v \in L^2(\Omega_R) \ | \ \nabla \cdot v = 0 \ \ \text{in} \ \ \Omega_{R} \, , \ \ v \cdot \nu = 0 \ \ \text{on} \ \
\partial \Omega_{R} \} \, , \\[6pt]
& G_{2}(\Omega_R) := \{v\in L^2(\Omega_R) \ | \ \exists g\in H^1(\Omega_R) \, , \ v=\nabla g \ \ \text{in} \ \ \Omega_{R} \}\, ,
\end{aligned}
$$
see \cite[Section 1]{foiastemam78} and \cite[Chapter 1]{ladyzhenskaya1969mathematical}. Recall that a vector field $v \in L^2_{\rm div}(\Omega_R)$
has
a normal trace in $H^{-1/2}(\partial \Omega_{R})$ defined through the generalized Green formula
$$
\langle v \cdot \nu , \varphi \rangle_{\partial \Omega_{R}} = \int_{\Omega_{R}} \varphi(\nabla \cdot v) + \int_{\Omega_{R}} v \cdot \nabla \varphi
\qquad \forall \varphi \in H^{1}(\Omega_{R}) \, ,
$$
where the ``boundary term'' $\langle\cdot,\cdot\rangle_{\partial \Omega_{R}}$ denotes the duality product between $H^{-1/2}(\partial \Omega_{R})$
and $H^{1/2}(\partial \Omega_{R})$. Likewise, recall that the (weak)
operator $\text{curl} : L^2(\Omega_R) \longrightarrow H^{-1}(\Omega_{R})$ is defined as
$$
\langle \text{curl}(v) , \varphi \rangle_{\Omega_{R}} = - \int_{\Omega_{R}} v \cdot \text{curl}(\varphi) \qquad \forall v \in L^2(\Omega_R) \, , \ \forall \varphi \in H_{0}^{1}(\Omega_{R}) \, ,
$$
where $\langle\cdot,\cdot\rangle_{\Omega_{R}}$ denotes the duality product between $H^{-1}(\Omega_{R})$ and $H_{0}^{1}(\Omega_{R})$.
\par
\begin{figure}[H]
	\begin{center}
\vspace*{-3mm}
	\includegraphics[scale=0.75]{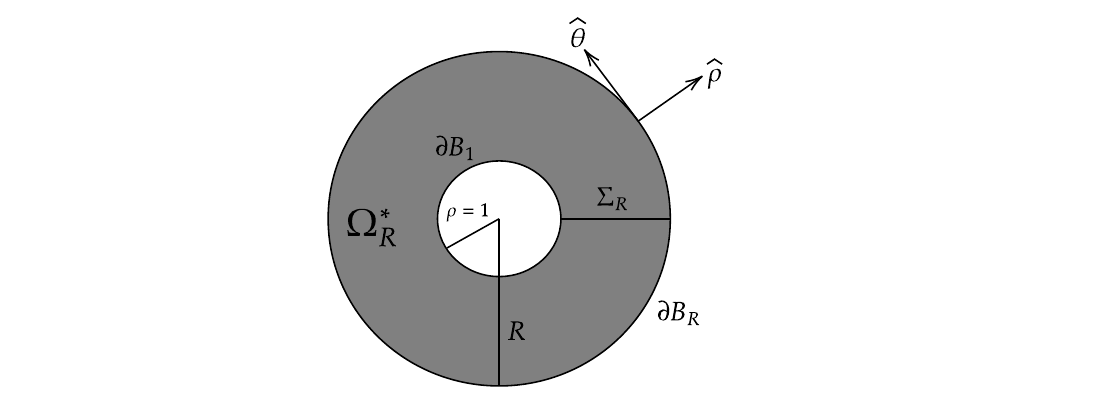}
	\end{center}
	\vspace*{-7mm}
	\caption{The cut annulus $\Omega^{*}_{R}$ in \eqref{cutannulus}.}\label{dom2}
\end{figure}

While $G_{2}(\Omega_R) \subset \text{Ker}(\text{curl})$, since the domain $\Omega_R$ is not simply connected, smooth irrotational vector fields are
not necessarily conservative (elements of $G_{2}(\Omega_R)$) and, therefore,
$$
G_{c}(\Omega_R) := G_{1}(\Omega_R) \cap \text{Ker}(\text{curl}) \neq\{0\} \, .
$$
The space $G_{c}(\Omega_R)$ can be characterized by choosing a smooth line $\Sigma_{R}\subset\Omega_R$ such that $\Omega_R\setminus
\overline{\Sigma_{R}}$ is simply connected. Roughly speaking,
$\Sigma_{R}$ connects $\partial B_1$ with $\partial B_R$ and cuts $\Omega_R$. We focus on the case where
$$
\Sigma_{R} := \{ \xi \in \Omega_{R} \ | \ \theta = 0 \} = \{ (x,y) \in \Omega_{R} \ | \ x > 0 \, , \ y=0\}
$$
and we explain below how to proceed for different choices of $\Sigma_{R}$, see Remark \ref{alternative} at the end of this section.
The open simply connected domain
\begin{equation} \label{cutannulus}
\Omega_R^* := \Omega_R\setminus \overline{\Sigma_R}
\end{equation}
has a {\it double boundary} at $\Sigma_R$: an upper boundary $\Sigma_R^+$ (when approached from the upper half plane $y>0$) and a lower boundary
$\Sigma_R^-$
(when approached from the lower half plane $y<0$), see Figure \ref{dom2}.
\noindent
For a continuous scalar function $q : \Omega_R^* \longrightarrow \mathbb{R}$, we denote by $[q]_{R} : \Sigma_R \longrightarrow \mathbb{R}$ the
\textit{jump function} when going from $\Sigma_R^-$ to $\Sigma_R^+$ (upwards). Precisely,
$$
[q]_{R}(x) := \lim_{y \to 0^{-}} q(x,y) - \lim_{y \to 0^{+}} q(x,y) \qquad \forall x \in (1,R) \, .
$$
In particular, if $[q]_{R} \equiv 0$, then $q$ can be extended continuously to the whole annulus $\Omega_R$. The following statement is a
consequence of some results by Foias \& Temam \cite{foiastemam78}:

\begin{proposition}\label{foiastemam}
The space $G_{c}(\Omega_R)$ has dimension $1$ and
$$
G_{c}(\Omega_R) = \textup{span} \left\{ \dfrac{1}{\rho} \, \widehat{\theta} \right\} \, .
$$
The orthogonal complement of $G_{c}(\Omega_R)$ within $G_{1}(\Omega_R)$ is characterized by
\begin{equation} \label{perpg}
G_{c}(\Omega_R)^{\perp} = \{ v \in G_{1}(\Omega_R) \ | \ \langle v \cdot \nu , 1 \rangle_{\Sigma_{R}} = 0 \} \, ,
\end{equation}
where $\langle\cdot,\cdot\rangle_{\Sigma_{R}}$ denotes the duality product between $H^{-1/2}(\Sigma_{R})$ and $H^{1/2}(\Sigma_{R})$.
\end{proposition}
\noindent
\begin{proof} From \cite[Lemmas 1.1-1.2-1.3]{foiastemam78} we know that the space $G_{c}(\Omega_R)$ has dimension $1$, and all the elements of
$G_{c}(\Omega_R)$ are proportional to the gradient of the unique (up to the addition of a constant) analytic function $\overline{q} \in
\mathcal{C}^{\infty}(\Omega_R^*) \cap L^{\infty}(\Omega_R)$ such that
\begin{equation}\label{firstchar}
\Delta \overline{q}=0 \ \mbox{ in } \ \Omega_R^*\, ,\qquad
\frac{\partial \overline{q}}{\partial\nu}=0 \ \mbox{ on } \ \partial\Omega_R\, ,\qquad [\overline{q}]_{R} = 1\, ,\qquad \left[\frac{\partial
\overline{q}}{\partial y}\right]_{R} = 0 \, .
\end{equation}
In our framework, a solution to \eqref{firstchar} is explicitly given by
$$
\overline{q}(x,y) = \dfrac{\theta}{2 \pi} = \left\{
\begin{array}{ll}
\displaystyle\frac{1}{2\pi}\arctan \left( \dfrac{y}{x} \right) \quad & \text{if} \ \ x>0 \ \ \text{and} \ \ y>0   \, , \\[14pt]
\displaystyle\frac{1}{2\pi}\arctan \left( \dfrac{y}{x} \right) + \frac12 \quad & \text{if} \ \ x < 0  \, , \\[14pt]
\displaystyle\frac{1}{2\pi}\arctan \left( \dfrac{y}{x} \right) + 1 \quad & \text{if} \ \ x>0 \ \ \text{and} \ \ y<0 \, ,
\end{array}
\right.
\qquad \forall (x,y) \in \Omega_{R} \, ,
$$
extended by continuity on the two segments where $x=0$, so that
$$\nabla \overline{q}(x,y)= \frac{1}{2 \pi \sqrt{x^2 + y^2}} \, \widehat{\theta}(x,y) \qquad \forall (x,y) \in \Omega_{R}^*\, .$$
From \cite[Lemma 1.4]{foiastemam78} we deduce that the orthogonal complement of $G_{c}(\Omega_R)$ within $G_{1}(\Omega_R)$ is the space given in
\eqref{perpg}.
\end{proof}

In Section \ref{uniquetc}, see \eqref{nstokesdebil}, we considered weak solutions to \eqref{nsstokes0}, which become strong solutions because
$f\in L^2(\Omega_R)$. Here we define a new class of solutions for which we need to introduce the following spaces of vector fields:
$$
\begin{aligned}
& \mathbb{H}(\Omega_{R}) := \left\{ v \in H^{1}_{0}(\Omega_{R}) \ \bigg| \ \int_{\Sigma_{R}} v \cdot \nu = 0 \right\} \, , \\[6pt]
& \mathbb{H}_{\sigma}(\Omega_{R}) := G_{c}(\Omega_R)^{\perp} \cap H^{1}_{0}(\Omega_{R}) = \left\{ v \in H^{1}_{0}(\Omega_{R}) \ \bigg|
\ \nabla \cdot v = 0 \ \ \text{in} \ \ \Omega_{R} \, , \ \int_{\Sigma_{R}} v \cdot \nu = 0 \right\} \, ,
\end{aligned}
$$
which, owing to the continuity of the trace operator, are closed subspaces of $H^{1}_{0}(\Omega_{R})$.

\begin{definition} \label{incompletedef}
Given $\omega\ge0$ and $f \in L^2(\Omega_R)$, we say that a vector field $u \in H^{1}_{*,\sigma}(\Omega_{R})$ is an
\textbf{incomplete solution} to \eqref{nsstokes0} if it verifies \eqref{nsstokes0}$_2$ in the trace sense and
\begin{equation} \label{nstokesdebilstar}
\int_{\Omega_{R}} \nabla u \cdot \nabla \varphi + \int_{\Omega_{R}} (u \cdot \nabla) u \cdot \varphi  = \int_{\Omega_{R}} f \cdot \varphi \qquad
\forall \varphi \in \mathbb{H}_{\sigma}(\Omega_{R})\, .
\end{equation}
\end{definition}
Since $\mathbb{H}_{\sigma}(\Omega_{R})\subset H^{1}_{0,\sigma}(\Omega_{R})$, every weak (strong) solution to \eqref{nsstokes0}
is also an incomplete
solution but, contrary to the standard framework \eqref{nstokesdebil}, no general regularity theory holds in the class of incomplete solutions, since Proposition \ref{foiastemam} shows that the space of test functions $\mathbb{H}_{\sigma}(\Omega_{R})$ lacks one dimension. In Theorem \ref{nonuniquetheo} below, we explicitly build infinitely many incomplete solutions (which are \textit{not} weak solutions) to the
unforced equation \eqref{nsstokes0} where the velocity field is smooth in $\Omega_{R}$ but the associated scalar pressure is {\it discontinuous}
along $\Sigma_{R}$, although its gradient has no jump.\par

Consider the space of scalar functions with zero mean value in $\Omega_{R}$, defined as
\begin{equation} \label{l02}
	L^2_0(\Omega_{R}) := \left\{ g\in L^2(\Omega_{R}) \ \Big | \ \int_{\Omega_{R}} g = 0 \right\} \, .
\end{equation}
For every vector field $\varphi \in \mathbb{H}(\Omega_{R})$ one certainly has $\nabla \cdot \varphi \in L^2_0(\Omega_{R})$, due to the Divergence Theorem. Reciprocally, in the spirit of \cite[Lemma 1]{pileckas1983spaces} (see also \cite[Lemma 1.8.13]{korobkov2024steady}), we have the following result concerning the inverse of the divergence operator:

\begin{lemma} \label{newbogovskii}
Given any scalar function $q \in L^2_0(\Omega_{R})$, there exists a vector field $J \in \mathbb{H}(\Omega_{R})$ such that
\begin{equation} \label{newbogeq}
 \nabla \cdot J = q \ \ \mbox{in} \ \ \Omega_{R} \qquad \text{and} \qquad \| \nabla J \|_{L^{2}(\Omega_{R})} \leq C_{*} \| q \|_{L^{2}(\Omega_{R})} \, ,
\end{equation}
	for some constant $C_{*} > 0$ that depends exclusively on $\Omega_{R}$.
\end{lemma}
\noindent
\begin{proof}
In what follows, $C > 0$ will always denote a generic constant that depends only on $\Omega_{R}$, but that may change from line to line.

Given any scalar function $q \in L^2_0(\Omega_{R})$, standard results such as \cite[Lemma 1.8.13]{korobkov2024steady} ensure the existence of a vector field $X \in H_{0}^{1}(\Omega_{R})$ such that
\begin{equation} \label{newbogeq2}
	\nabla \cdot X = q \ \ \mbox{in} \ \ \Omega_{R} \qquad \text{and} \qquad \| \nabla X \|_{L^{2}(\Omega_{R})} \leq C \| q \|_{L^{2}(\Omega_{R})} \, ,
\end{equation}
which can be written as
$$
X(\xi) = X_{\rho}(\rho, \theta) \widehat{\rho} + X_{\theta}(\rho,\theta) \widehat{\theta} \qquad \text{for a.e. }\xi \in \Omega_{R}
\, ,
$$
for some scalar functions $X_{\rho}, X_{\theta}  \in H_{0}^{1}(\Omega_{R})$. Now, consider the function
\begin{equation} \label{fung}
g(\rho) := \dfrac{6}{(R - 1)^{3}} (\rho - 1)(R - \rho) \quad \forall \rho \in [1,R] \quad \Longrightarrow \quad g \in H_{0}^{1}(1,R) \cap \mathcal{C}^{\infty}([1,R]) \, , \ \ \int_{1}^{R} g(\rho) \, d\rho = 1 \, .
\end{equation}
Then, define the vector field $J \in H_{0}^{1}(\Omega_{R})$ by
$$
J(\xi) := X_{\rho}(\rho, \theta) \widehat{\rho} + \left[ X_{\theta}(\rho, \theta) - \left( \int_{\Sigma_{R}} X \cdot \nu \right) g(\rho) \right] \widehat{\theta} \qquad \text{for a.e. }\xi \in \Omega_{R} \, ,
$$
which, due to \eqref{newbogeq2}$_1$-\eqref{fung}, is an element of $\mathbb{H}(\Omega_{R})$ such that $\nabla \cdot J = \nabla \cdot X = q$ in $\Omega_{R}$. Moreover, since the trace inequality and \eqref{newbogeq2}$_2$ entail
$$
\left| \int_{\Sigma_{R}} X \cdot \nu \right| \leq C \| X \|_{L^{2}(\Sigma_{R})} \leq C \| \nabla X \|_{L^{2}(\Omega_{R})} \leq C \| q \|_{L^{2}(\Omega_{R})} \, ,
$$
we conclude that $J \in \mathbb{H}(\Omega_{R})$ verifies all the properties stated in \eqref{newbogeq}.
\end{proof}

We are now in position to explain how a pressure $p\in L^2(\Omega_R)/\mathbb{R}$, associated to an incomplete solution to \eqref{nsstokes0}, arises.

\begin{theorem} \label{incompletepressure}
For $\omega\ge0$ and $f \in L^2(\Omega_R)$, let $u \in H^{1}_{*,\sigma}(\Omega_{R})$ be an incomplete solution to \eqref{nsstokes0}.
Then, there exists a unique scalar pressure $p \in L^{2}(\Omega_{R})/\mathbb{R}$ such that
\begin{equation} \label{nstokesdebilstarpres}
\int_{\Omega_{R}} \nabla u \cdot \nabla \varphi + \int_{\Omega_{R}} (u \cdot \nabla) u \cdot \varphi  - \int_{\Omega_{R}} p (\nabla \cdot \varphi)
= \int_{\Omega_{R}} f \cdot \varphi \qquad \forall \varphi \in \mathbb{H}(\Omega_{R}) \, .
\end{equation}
\end{theorem}
\noindent
\begin{proof} We follow closely the proof of \cite[Theorem III.5.3]{galdi2011introduction} (see also \cite[Lemma 2.1]{girault2012finite}).
	
Since for every vector field $\varphi \in \mathbb{H}(\Omega_{R})$ we have  $\nabla \cdot \varphi \in L^2_0(\Omega_{R})$, it suffices to prove the existence of a unique $p \in L_{0}^{2}(\Omega_{R})$ satisfying \eqref{nstokesdebilstarpres}.
Denote
by $\mathbb{H}^{-1}(\Omega_{R})$ the dual space of $\mathbb{H}(\Omega_{R})$ and by $\langle\cdot,\cdot\rangle_{\Omega_{R}}^*$ the
duality product
between $\mathbb{H}^{-1}(\Omega_{R})$ and $\mathbb{H}(\Omega_{R})$.

A simple density argument, based on integration by parts (see, for example, \cite[Lemma 2.1]{gazgruspe2026}), furnishes the inequality
\begin{equation} \label{grunau}
\| \nabla \cdot \varphi \|_{L^{2}(\Omega_{R})} \leq \| \nabla \varphi \|_{L^{2}(\Omega_{R})} \qquad \forall \varphi \in H_{0}^{1}(\Omega_{R}) \, .
\end{equation}
Then, the (weak) gradient
$-\nabla(\cdot): L^2_0(\Omega_{R}) \longrightarrow \mathbb{H}^{-1}(\Omega_{R})$ is the bounded linear application defined by
$$
\langle -\nabla q , \varphi \rangle_{\Omega_{R}}^* := \int_{\Omega_{R}} q (\nabla \cdot \varphi) \qquad \forall q \in L^2_0(\Omega_{R}) \, , \ \ \forall \varphi \in \mathbb{H}(\Omega_{R}) \, ,
$$
so that it is the adjoint of the (strong) divergence operator $\text{div} : \mathbb{H}(\Omega_{R})\longrightarrow L^2_0(\Omega_{R})$. From \eqref{grunau} we deduce the following inequalities:
$$
\| \nabla q \|_{\mathbb{H}^{-1}(\Omega_{R})} \leq \| q \|_{L^{2}(\Omega_{R})} \quad \forall q \in L^2_0(\Omega_{R}) \, ; \qquad \| \nabla \cdot \varphi \|_{L^{2}(\Omega_{R})} \leq \| \nabla \varphi \|_{L^{2}(\Omega_{R})} \qquad \forall \varphi \in \mathbb{H}(\Omega_{R}) \, .
$$
Moreover, Lemma \ref{newbogovskii} entails that such divergence operator is surjective from $\mathbb{H}(\Omega_{R})$ to $L^2_0(\Omega_{R})$ (thus, its range is closed in $L^2_0(\Omega_{R})$). Therefore, the Closed Range Theorem \cite[Chapter VII]{yosida2012functional} yields
$$
\text{Range}(-\nabla) = \text{Ker}(\text{div})^{\perp} = \mathbb{H}_{\sigma}(\Omega_{R})^{\perp} := \{ \mathcal{Q} \in
\mathbb{H}^{-1}(\Omega_{R}) \ | \ \langle \mathcal{Q} , \varphi \rangle_{\Omega_{R}}^* = 0 \quad \forall \varphi \in
\mathbb{H}_{\sigma}(\Omega_{R}) \} \, .
$$
Given an incomplete solution $u \in H^{1}_{*,\sigma}(\Omega_{R})$ of \eqref{nsstokes0}, define the linear functional
$\mathcal{F}_u : \mathbb{H}(\Omega_{R}) \longrightarrow \mathbb{R}$ by
$$
\mathcal{F}_u(\varphi) := \int_{\Omega_{R}} \nabla u \cdot \nabla \varphi + \int_{\Omega_{R}} (u \cdot \nabla) u \cdot \varphi -
\int_{\Omega_{R}} f \cdot \varphi \qquad \forall \varphi \in \mathbb{H}(\Omega_{R}) \, .
$$
Applying the H\"older, Poincaré and Sobolev inequalities in $\Omega_{R}$ (see \eqref{sobolevconstants11}) we infer
$$
| \mathcal{F}_u(\varphi) | \leq \left( \dfrac{1}{\mathcal{S}_{4}} \| \nabla u \|^{2}_{L^{2}(\Omega_{R})} +  \| \nabla u \|_{L^{2}(\Omega_{R})} + \dfrac{1}{\sqrt{\mathcal{S}_{2}}} \| f \|_{L^{2}(\Omega_{R})} \right) \| \nabla \varphi \|_{L^{2}(\Omega_{R})} \qquad \forall \varphi \in \mathbb{H}(\Omega_{R}) \, ,
$$
so that $\mathcal{F}_u \in \mathbb{H}^{-1}(\Omega_{R})$. Moreover, owing  to \eqref{nstokesdebilstar}, we deduce that $\mathcal{F}_u\in
\text{Range}(-\nabla)$. This ensures the existence of a scalar pressure $p \in L_{0}^{2}(\Omega_{R})$ verifying \eqref{nstokesdebilstarpres}.
\par
Suppose there exists another scalar function $q \in L_{0}^{2}(\Omega_{R})$ verifying \eqref{nstokesdebilstarpres}, thereby implying
\begin{equation} \label{bog1}
\int_{\Omega_{R}} (p-q) (\nabla \cdot \varphi) = 0 \qquad \forall \varphi \in \mathbb{H}(\Omega_{R}) \, .
\end{equation}
As $p - q \in L_{0}^{2}(\Omega_{R})$, Lemma \ref{newbogovskii} ensures the existence of a vector field $J \in \mathbb{H}(\Omega_{R})$ such that $\nabla \cdot J = p-q$ in $\Omega_{R}$. Taking $\varphi = J$ in \eqref{bog1} gives us $p \equiv q$ almost everywhere in $\Omega_{R}$, which concludes the proof.
\end{proof}

Let $v \in H^{1}(\Omega_{R})$ be a vector field, decomposed as
$$
v(\xi) = v_{\rho}(\rho, \theta) \widehat{\rho} + v_{\theta}(\rho,\theta) \widehat{\theta} \qquad \text{for a.e. }\xi \in \Omega_{R}
\, ,
$$
for some scalar functions $v_{\rho}, v_{\theta}  \in H^{1}(\Omega_{R})$. We define its \textit{flux across $\Sigma_{R}$} as the quantity
$$
\Phi_{v} := \int_{\Sigma_{R}} v \cdot \nu = \int_{1}^{R} v_{\theta}(\rho,0) \, d\rho \, .
$$
Given $\omega \geq 0$, a vector field $v \in H^{1}(\Omega_{R})$ satisfying
\begin{equation} \label{bog7}
	\nabla\cdot v=0 \ \text{ in } \ \Omega_{R} \, , \qquad v = (0,0) \ \text{ on } \ \partial B_{1} \, , \qquad v = \omega \, \widehat{\theta} \
\text{ on } \ \partial B_{R}\, ,
\end{equation}
has \textit{constant flux} across each angular sector of $\Omega_{R}$. Indeed, for any $\alpha \in [0, 2 \pi)$, define the segment
$$
\Sigma^{\alpha}_{R} := \{ \xi \in \overline{\Omega_{R}} \ | \ 1 \leq \rho \leq R \, , \ \theta = \alpha \} \, .
$$
By the Divergence Theorem and the boundary conditions in \eqref{bog7} we have that
\begin{equation} \label{sameflux}
\int_{\Sigma_{R}} v \cdot \nu = \int_{\Sigma^{\alpha}_{R}} v \cdot \nu \qquad \forall \alpha \in [0, 2 \pi) \, .
\end{equation}

The main result of this section, providing infinitely many incomplete solutions to \eqref{nsstokes0}, is obtained by
showing the existence of an incomplete solution to \eqref{nsstokes0} with a \textit{prescribed} flux $\Phi$, for every $\Phi \in
\mathbb{R}$.

\begin{theorem} \label{nonuniquetheomain}
For any $\omega \geq 0$ and $f \in \mathcal{R}[L^{2}(\Omega_{R})]$, there exist infinitely many rotationally-invariant incomplete solutions to
\eqref{nsstokes0}.
\end{theorem}
\noindent
\begin{proof} We first build a \textit{flux carrier} as an incomplete solution of a (linear) Stokes problem.
For $\omega \geq 0$ and $\Phi \in \mathbb{R}$, put
\begin{equation} \label{lambdastar}
	\lambda_{*} := \dfrac{8\Phi (R^2 - 1) - 4 \omega R \left( R^2 - 2\log(R) - 1 \right)}{(R^2 - 1)^2 - 4 R^2 (\log(R))^{2}} \, .
\end{equation}
By direct computation one can check that the couple
\begin{equation}\label{fluxcarrier}
V_{*}(\xi) := \left[ \dfrac{R}{R^2 - 1} \left( \omega + \dfrac{\lambda_{*} R}{2} \log(R) \right)  \left( \rho - \dfrac{1}{\rho} \right) -
\dfrac{\lambda_{*} \rho}{2} \log(\rho) \right] \widehat{\theta} \,, \qquad  Q_{*}(\xi) := - \lambda_{*}\theta \qquad \forall \xi \in
\overline{\Omega_{R}}
\end{equation}
is the unique incomplete solution (up to an additive constant for $Q_{*}$) of the unforced Stokes system
\begin{equation}\label{stokesVinc}
	\left\{
	\begin{aligned}
		& -\Delta v+\nabla q=0 \, , \ \quad  \nabla\cdot v=0 \ \ \mbox{ in } \ \ \Omega_{R} \, , \\[3pt]
		& v = \omega \, \widehat{\theta} \ \text{ on } \ \partial B_R \, , \qquad v = (0,0) \ \text{ on } \ \partial B_{1} \, ,
	\end{aligned}
	\right.
\end{equation}
having flux $\Phi$. Indeed,
\begin{equation}\label{lambdastarwf}
\int_{\Omega_R}\nabla V_{*}\cdot\nabla\varphi-\int_{\Omega_{R}}Q_{*}(\nabla\cdot\varphi)=0\qquad\forall\varphi\in
\mathbb{H}(\Omega_R)
\end{equation}
and, by linearity, the unique solvability of \eqref{stokesVinc} in the class of incomplete solutions having flux $\Phi$ follows,
since the difference of any two such solutions is an element of $\mathbb{H}(\Omega_{R})$.\par
For any $\Phi\in\mathbb{R}$, let $V_{*}\in\mathcal{C}^{\infty}(\overline{\Omega_{R}})$ be as in \eqref{fluxcarrier}. We now
prove the existence of a rotationally-invariant incomplete solution $u\in\mathcal{R}[H^{1}_{*,\sigma}(\Omega_{R})]$ to
\eqref{nsstokes0} having flux $\Phi$. This amounts to showing the existence of
$\widehat{u}\in\mathcal{R}[\mathbb{H}_{\sigma}(\Omega_{R})]$ such that
\begin{equation}\label{oseenquasi}
	\begin{aligned}
		& \int_{\Omega_{R}} \nabla \widehat{u} \cdot \nabla \varphi + \int_{\Omega_{R}} (\widehat{u} \cdot \nabla)\widehat{u} \cdot \varphi +
\int_{\Omega_{R}} (\widehat{u} \cdot \nabla)V_{*} \cdot \varphi + \int_{\Omega_{R}} (V_{*} \cdot \nabla)\widehat{u} \cdot \varphi + \int_{\Omega_{R}} (V_{*} \cdot \nabla)V_{*} \cdot \varphi \\[6pt]
		& \hspace{-4mm} = \int_{\Omega_{R}} f \cdot \varphi \qquad \forall \varphi \in
\mathbb{H}_{\sigma}(\Omega_{R}) \, ,
	\end{aligned}
\end{equation}
and the solution will then be $u=\widehat{u} + V_{*}$, see again \eqref{lambdastarwf}.

We start by showing the existence of a unique vector field $\widehat{u} \in \mathbb{H}_{\sigma}(\Omega_{R})$ satisfying
\begin{equation}\label{oseenquasi2}
\int_{\Omega_{R}} \nabla \widehat{u} \cdot \nabla \varphi + \int_{\Omega_{R}} \left[ (\widehat{u} \cdot \nabla)V_{*} + (V_{*} \cdot \nabla)\widehat{u} \right] \cdot \varphi  = \int_{\Omega_{R}} \left[ f - (V_{*} \cdot \nabla)V_{*} \right] \cdot \varphi \qquad \forall \varphi \in \mathbb{H}_{\sigma}(\Omega_{R}) \, .
\end{equation}
Indeed, from the H\"older, Poincaré and Sobolev inequalities in $\Omega_{R}$ (see \eqref{sobolevconstants11}) it follows that
$$
\left\{
\begin{aligned}
& (v, \varphi) \in \mathbb{H}_{\sigma}(\Omega_{R}) \times \mathbb{H}_{\sigma}(\Omega_{R}) \longmapsto \int_{\Omega_{R}} \nabla v \cdot \nabla \varphi + \int_{\Omega_{R}} \left[ (v \cdot \nabla)V_{*} + (V_{*} \cdot \nabla)v \right] \cdot \varphi \, ,\\[6pt]
& \varphi \in \mathbb{H}_{\sigma}(\Omega_{R}) \longmapsto \int_{\Omega_{R}} \left[ f - (V_{*} \cdot \nabla)V_{*} \right] \cdot \varphi \, ,
\end{aligned}
\right.
$$
constitute, respectively, continuous bilinear and linear applications on $\mathbb{H}_{\sigma}(\Omega_{R}) \times \mathbb{H}_{\sigma}(\Omega_{R})$ and $\mathbb{H}_{\sigma}(\Omega_{R})$. Therefore, the Lax-Milgram Theorem guarantees the existence of a unique vector field $\widehat{u} \in \mathbb{H}_{\sigma}(\Omega_{R})$ satisfying \eqref{oseenquasi2}, which we decompose as
\begin{equation} \label{rotinvariant2hat}
\widehat{u}(\xi) = u_{\rho}(\rho, \theta) \widehat{\rho} + u_{\theta}(\rho,\theta) \widehat{\theta} \qquad \text{for a.e. }\xi \in \Omega_{R}
\, ,
\end{equation}
for some scalar functions $u_{\rho}, u_{\theta}  \in H_{0}^{1}(\Omega_{R})$. Given any angle $\phi \in \mathbb{R}$ define, as in \eqref{rotinvariant}, the vector field
\begin{equation} \label{rotinvariant2}
	\widehat{u}_{\phi}(\xi) := \mathcal{R}(\phi)^{\top} \widehat{u}(\mathcal{R}(\phi) \xi) = u_{\rho}(\rho, \theta + \phi) \widehat{\rho}(\theta) + u_{\theta}(\rho, \theta + \phi) \widehat{\theta}(\theta)  \qquad
	\forall \xi \in \Omega_{R} \, ,
\end{equation}
so that $\widehat{u}_{\phi} \in H^{1}_{0,\sigma}(\Omega_{R})$. Moreover, recalling \eqref{sameflux}, we have
$$
\int_{\Sigma_{R}} \widehat{u}_{\phi} \cdot \nu = \int_{\Sigma^{\phi}_{R}} \widehat{u} \cdot \nu = \int_{\Sigma_{R}} \widehat{u} \cdot \nu = 0 \, ,
$$
thus implying $\widehat{u}_{\phi} \in \mathbb{H}_{\sigma}(\Omega_{R})$. Similarly, fixing $\varphi \in \mathbb{H}_{\sigma}(\Omega_{R})$, there holds $\varphi^{(-\phi)} \in \mathbb{H}_{\sigma}(\Omega_{R})$, with
$$
\varphi^{(-\phi)}(\xi) := \mathcal{R}(-\phi)^{\top} \varphi(\mathcal{R}(-\phi) \xi) = \varphi_{\rho}(\rho, \theta - \phi) \widehat{\rho}(\theta) + \varphi_{\theta}(\rho, \theta - \phi) \widehat{\theta}(\theta) \qquad
\forall \xi \in \Omega_{R} \, ,
$$
for some scalar functions $\varphi_{\rho}, \varphi_{\theta}  \in H_{0}^{1}(\Omega_{R})$. Therefore, the integral identity \eqref{oseenquasi2} yields
\begin{equation} \label{rotinvariant3}
\int_{\Omega_{R}} \nabla \widehat{u} \cdot \nabla \varphi^{(-\phi)} + \int_{\Omega_{R}} \left[ (\widehat{u} \cdot \nabla)V_{*} + (V_{*} \cdot \nabla)\widehat{u} \right] \cdot \varphi^{(-\phi)}  = \int_{\Omega_{R}} \left[ f - (V_{*} \cdot \nabla)V_{*} \right] \cdot \varphi^{(-\phi)} \, .
\end{equation}
Since $\widehat{u}, \varphi \in H^{1}(\Omega_{R})$, they are $2\pi$-periodic in the variable $\theta \in [0, 2 \pi]$, together with their gradients. As $ V_{*} \in \mathcal{R}[H^{1}(\Omega_{R})]$ and $f \in \mathcal{R}[L^{2}(\Omega_{R})]$, we deduce, from \eqref{rotinvariant3} and a simple change of variables,
\begin{equation} \label{rotinvariant4}
	\int_{\Omega_{R}} \nabla \widehat{u}_{\phi} \cdot \nabla \varphi + \int_{\Omega_{R}} \left[ (\widehat{u}_{\phi} \cdot \nabla)V_{*} + (V_{*} \cdot \nabla)\widehat{u}_{\phi} \right] \cdot \varphi  = \int_{\Omega_{R}} \left[ f - (V_{*} \cdot \nabla)V_{*} \right] \cdot \varphi \qquad \forall \varphi \in \mathbb{H}_{\sigma}(\Omega_{R}) \, .
\end{equation}
By uniqueness, there holds $\widehat{u} \equiv \widehat{u}_{\phi}$ in $\Omega_{R}$ for every $\phi \in \mathbb{R}$, so that $\widehat{u}\in\mathcal{R}[\mathbb{H}_{\sigma}(\Omega_{R})]$. Thus, since $u_{\theta}$ is independent
of $\theta$ (see again \eqref{rotinvariant2hat}), the incompressibility condition becomes
$$
\dfrac{d}{d\rho}(\rho \, u_{\rho}(\rho)) = 0 \quad \Longrightarrow \quad u_{\rho}(\rho) = \dfrac{A}{\rho} \qquad\text{for a.e.
}\rho\in(1,R) \, ,
$$
for some constant $A \in \mathbb{R}$. Since $\widehat{u}$ vanishes on $\partial \Omega_{R}$, we must have $A=0$, yielding $\widehat{u}(\xi) = u_{\theta}(\rho) \widehat{\theta}$ for a.e. $\xi \in \Omega_{R}$. Subsequently, defining the scalar function $U_{*} \in H^{1}(\Omega_{R})$ by
\begin{equation} \label{rotinvariant5}
U_{*}(\xi) := - \int_{1}^{\rho} \dfrac{1}{t} (u_{\theta}(t))^{2} \, dt \qquad \forall \xi \in \Omega_{R} \, ,
\end{equation}
we certainly have
\begin{equation} \label{rotinvariant6}
((\widehat{u} \cdot \nabla) \widehat{u})(\xi) = - \dfrac{1}{\rho} (u_{\theta}(\rho))^{2} \widehat{\rho} = \nabla U_{*}(\xi) \qquad \forall \xi \in \Omega_{R} \, .
\end{equation}
On the other hand, exactly as in the proof of Theorem \ref{incompletepressure} we deduce the existence of a unique scalar pressure $\widehat{p} \in L_{0}^{2}(\Omega_{R})$ such that
\begin{equation}\label{oseenquasi2pres}
	\int_{\Omega_{R}} \nabla \widehat{u} \cdot \nabla \varphi + \int_{\Omega_{R}} \left[ (\widehat{u} \cdot \nabla)V_{*} + (V_{*} \cdot \nabla)\widehat{u} \right] \cdot \varphi - \int_{\Omega_{R}} \widehat{p}(\nabla \cdot \varphi) = \int_{\Omega_{R}} \left[ f - (V_{*} \cdot \nabla)V_{*} \right] \cdot \varphi \quad \forall \varphi \in \mathbb{H}(\Omega_{R}) \, ,
\end{equation}
which, owing to \eqref{rotinvariant5}-\eqref{rotinvariant6} and an integration by parts, is equivalent to
\begin{equation}\label{oseenquasi2pres2}
	\begin{aligned}
		& \int_{\Omega_{R}} \nabla \widehat{u} \cdot \nabla \varphi + \int_{\Omega_{R}} \left[ (\widehat{u} \cdot \nabla)\widehat{u} + (\widehat{u} \cdot \nabla)V_{*} + (V_{*} \cdot \nabla)\widehat{u} \right] \cdot \varphi - \int_{\Omega_{R}} (\widehat{p} - U_{*})(\nabla \cdot \varphi) \\[6pt]
		& \hspace{-4mm} = \int_{\Omega_{R}} \left[ f - (V_{*} \cdot \nabla)V_{*} \right] \cdot \varphi \qquad \forall \varphi \in
		\mathbb{H}(\Omega_{R}) \, .
	\end{aligned}
\end{equation}
In turn, \eqref{oseenquasi2pres2} ensures the validity of \eqref{oseenquasi} by considering test functions in $\mathbb{H}_{\sigma}(\Omega_{R})$.

In conclusion, for every $\Phi \in \mathbb{R}$ there exists (at least) one rotationally-invariant
incomplete solution to \eqref{nsstokes0} having flux $\Phi$.
\end{proof}

\begin{remark}
Theorem \ref{nonuniquetheomain} explains the uniqueness of strong rotationally-invariant solutions proved in Theorem \ref{UVconstant0},
regardless of the magnitude of the angular velocity $\omega$. Recalling the bifurcation result in the region between two concentric
unbounded 3D cylinders \cite{velte1966stabilitat}, having $\Omega_R$ as planar cross-section, Theorem \ref{UVconstant0} states that if the pressure
is not allowed to vary along the vertical axis, then the rotationally-invariant solution cannot bifurcate. On the other hand, in view of Proposition
\ref{foiastemam}, Theorem \ref{nonuniquetheomain} suggests that, perhaps, an {\em helicoidal} solution in the region between the two cylinders
may be constructed through the Riemann surfaces of the multivalued function $\xi \mapsto \theta$, and possibly related to the leapfrogging vortex ring phenomenon \cite{davila2024leapfrogging}.
\end{remark}

\begin{remark}
It is natural to wonder whether one can relax the assumption of rotational invariance in Theorem \ref{nonuniquetheomain}, that is, to prove existence of infinitely many (not necessarily rotationally-invariant) incomplete solutions to \eqref{nsstokes0}. Notice that the involved boundary data \eqref{nsstokes0}$_2$ verify the \textbf{stringent outflow condition}, as described in \cite[Section 3.3.1]{korobkov2024steady}. Then, following standard methods in the theory of the stationary Navier-Stokes equations under non-homogeneous Dirichlet boundary conditions \cite[Chapter 5, Section 2]{ladyzhenskaya1969mathematical}, one is led to build a flux carrier that, apart from having flux $\Phi \in \mathbb{R}$ across $\Sigma_{R}$, must satisfy the Leray-Hopf inequality, see \cite[Lemma 3.3.1]{korobkov2024steady}. Such construction, albeit challenging and meaningful, is outside the scope of the present article.
\end{remark}

An interesting situation arises when
$$
f(\xi) = \dfrac{\lambda}{\rho} \, \widehat{\theta} \qquad \forall \xi \in \overline{\Omega_{R}} \, , \ \text{for some} \ \lambda \in \mathbb{R} \,
,
$$
and the generalized Taylor-Couette flow \eqref{gtc1} becomes
\begin{equation} \label{gtc2}
	\left\{
	\begin{aligned}
		& v_{*}(\xi) = \left[ \dfrac{R}{R^2 - 1} \left( \omega + \dfrac{\lambda R}{2} \log(R) \right)  \left( \rho - \dfrac{1}{\rho} \right) -
\dfrac{\lambda \rho}{2} \log(\rho) \right] \widehat{\theta} \\[6pt]
		& q_{*}(\xi) = \int_{1}^{\rho} \dfrac{1}{t} \left[ \dfrac{R}{R^2 - 1} \left( \omega + \dfrac{\lambda R}{2} \log(R) \right)  \left( t -
\dfrac{1}{t} \right) - \dfrac{\lambda t}{2} \log(t) \right]^{2} dt
	\end{aligned}
	\right. \qquad \forall \xi \in \overline{\Omega_{R}}\,.
\end{equation}
In the particular case when $f=0$ we can prove Theorem \ref{nonuniquetheomain} in a direct way, making use of the explicit generalized
Taylor-Couette flow \eqref{gtc2} with $\lambda=0$. This allows to show that a smooth rotationally-invariant incomplete solution to
\eqref{nsstokes0} may have an associated pressure that is not necessarily rotationally invariant, and not even continuous in $\Omega_{R}$.

\begin{theorem} \label{nonuniquetheo}
For any $\omega \geq 0$ there exist infinitely many incomplete solutions to the unforced system \eqref{nsstokes0} (with $f=0$),
given by the family
$$
\left\{
\begin{aligned}
	& v_{\lambda}(\xi) = \left[ \dfrac{R}{R^2 - 1} \left( \omega + \dfrac{\lambda R}{2} \log(R) \right)  \left( \rho - \dfrac{1}{\rho} \right) -
\dfrac{\lambda \rho}{2} \log(\rho) \right] \widehat{\theta} \\[6pt]
	& q_{\lambda}(\xi) = \int_{1}^{\rho} \dfrac{1}{t} \left[ \dfrac{R}{R^2 - 1} \left( \omega + \dfrac{\lambda R}{2} \log(R) \right)  \left( t -
\dfrac{1}{t} \right) - \dfrac{\lambda t}{2} \log(t) \right]^{2} dt - \lambda \theta
\end{aligned}
\right. \qquad \forall \xi \in \overline{\Omega_{R}} \, , \ \text{with} \ \lambda \in \mathbb{R} \, .
$$
\end{theorem}
\noindent
\begin{proof} Given $\varepsilon \in (0,\pi)$, consider the domain
	\begin{equation} \label{bog4}
		\Omega^{\varepsilon}_{R} := \{ \xi \in \Omega_{R} \ | \ \varepsilon < \theta < 2\pi - \varepsilon \} \, ,
	\end{equation}
	so that $\partial\Omega^\varepsilon_R=\partial B^\varepsilon_1\cup\partial
	B^\varepsilon_R\cup\Sigma^{\varepsilon,1}_{R}\cup\Sigma^{\varepsilon,2}_R$, where
	$$
	\begin{array}{ll}
		\partial B^{\varepsilon}_{1} := \{ \xi \in \mathbb{R}^{2} \ | \ \rho=1 \, , \quad \varepsilon < \theta < 2\pi - \varepsilon \} \, , \qquad
		& \partial B^{\varepsilon}_{R} := \{ \xi \in \mathbb{R}^{2} \ | \ \rho=R \, , \quad \varepsilon < \theta < 2\pi - \varepsilon \} \, , \\[6pt]
		\Sigma^{\varepsilon,1}_{R} := \{ \xi \in \mathbb{R}^{2} \ | \ \theta=\varepsilon \, , \quad 1 \leq \rho \leq R \} \, , \qquad &
		\Sigma^{\varepsilon,2}_{R} := \{ \xi \in \mathbb{R}^{2} \ | \ \theta=2\pi - \varepsilon \, , \quad 1 \leq \rho \leq R \} \, ,
	\end{array}
	$$
	see Figure \ref{pacman}.
	\begin{figure}[H]
		\begin{center}
			\vspace*{-3mm}
			\includegraphics[scale=0.75]{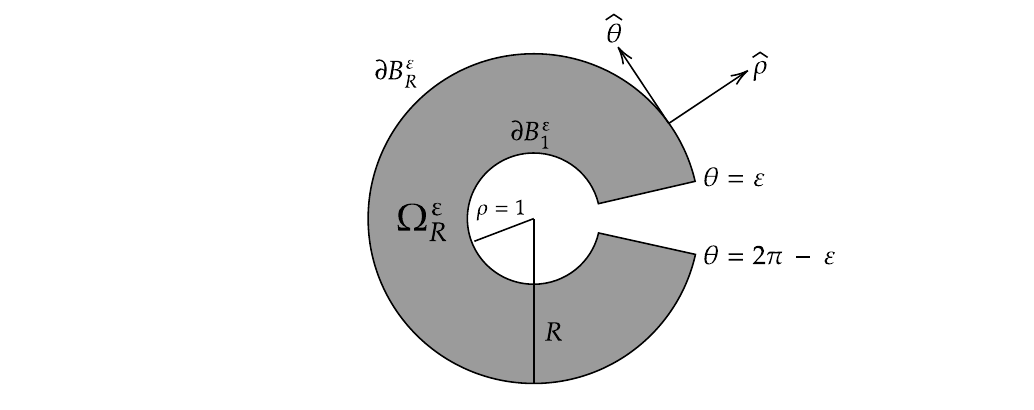}
		\end{center}
		\vspace*{-7mm}
		\caption{The domain $\Omega^{\varepsilon}_{R}$ in \eqref{bog4}.}\label{pacman}
	\end{figure}
For any vector field $\varphi \in H_{0}^{1}(\Omega_{R})$, an integration by parts yields
$$
\int_{\Omega^{\varepsilon}_{R}} \theta (\nabla \cdot \varphi) = - \int_{\Omega^{\varepsilon}_{R}} \dfrac{1}{\rho} \, \widehat{\theta} \cdot
\varphi + (2\pi - \varepsilon) \int_{\Sigma^{\varepsilon,2}_{R}} \varphi \cdot \widehat{\theta} - \varepsilon \int_{\Sigma^{\varepsilon,1}_{R}}
\varphi \cdot \widehat{\theta} \qquad \forall \varepsilon \in (0,\pi) \, .
$$
By the Dominated Convergence Theorem, if we take the limit as $\varepsilon \to 0^{+}$ in this identity, we obtain
$$
\int_{\Omega_{R}} \theta (\nabla \cdot \varphi) = 2\pi \int_{\Sigma_{R}} \varphi \cdot \widehat{\theta} - \int_{\Omega_{R}} \dfrac{1}{\rho} \,
\widehat{\theta} \cdot \varphi \qquad \forall \varphi \in H_{0}^{1}(\Omega_{R}) \, .
$$
In particular we have
\begin{equation} \label{dirac}
\nabla \theta = \dfrac{1}{\rho} \, \widehat{\theta} - 2 \pi \delta_{\Sigma_{R}} \widehat{e}_{2} \quad \text{in distributional sense}\, ,
\end{equation}
and also
\begin{equation} \label{bog5}
	\int_{\Omega_{R}} \theta (\nabla \cdot \varphi) = - \int_{\Omega_{R}} \dfrac{1}{\rho} \, \widehat{\theta} \cdot \varphi \qquad \forall \varphi
\in \mathbb{H}(\Omega_{R}) \, .
\end{equation}

Given $\omega \geq 0$ and $\lambda \in \mathbb{R}$, let $(v_{*},q_{*}) \in \mathcal{C}^\infty(\overline{\Omega_{R}}) \times
\mathcal{C}^\infty(\overline{\Omega_{R}})$, be the generalized Taylor-Couette flow defined in \eqref{gtc2}, which satisfies
$$
\left\{
\begin{aligned}
	& -\Delta v_{*}+ (v_{*}\cdot\nabla)v_{*}+ \nabla q_{*}=\dfrac{\lambda}{\rho}\, \widehat{\theta}\, ,\qquad \nabla\cdot v_{*}=0 \quad\mbox{in }
\Omega_{R}\, , \\[3pt]
	&  v_{*}= \omega \, \widehat{\theta} \quad\mbox{on }\partial B_R\, , \qquad v_{*} = (0,0) \quad\mbox{on }\partial B_{1} \, ,
\end{aligned}
\right.
$$
or equivalently,
\begin{equation} \label{bog6}
\int_{\Omega_{R}} \nabla v_{*} \cdot \nabla \varphi + \int_{\Omega_{R}} (v_{*} \cdot \nabla) v_{*} \cdot \varphi  - \int_{\Omega_{R}} q_{*}
(\nabla \cdot \varphi) = \int_{\Omega_{R}} \dfrac{\lambda}{\rho} \, \widehat{\theta} \cdot \varphi \qquad \forall \varphi \in
H_{0}^{1}(\Omega_{R}) \, .
\end{equation}
Identities \eqref{bog5}-\eqref{bog6} imply that
$$
	\int_{\Omega_{R}} \nabla v_{*} \cdot \nabla \varphi + \int_{\Omega_{R}} (v_{*} \cdot \nabla) v_{*} \cdot \varphi  - \int_{\Omega_{R}} (q_{*} -
\lambda \theta) (\nabla \cdot \varphi) = 0 \qquad \forall \varphi \in \mathbb{H}(\Omega_{R}) \, ,
$$
so that $v_{*}$ is also an incomplete solution to \eqref{nsstokes0} (with $f=0$) with associated pressure $q_{*}
- \lambda \theta$, which is a \textit{discontinuous} function on $\Sigma_{R}$. Hence, since $q_{*}$ is rotationally invariant, we infer
$q_{*} - \lambda \theta \notin H^{1}(\Omega_{R})$.
\end{proof}

\begin{remark}
As an example, given $\lambda \in \mathbb{R}$, the pair
$$
\left\{
\begin{aligned}
& v_{\lambda}(\xi) = \lambda \left[ \dfrac{R^{2} \log(R)}{2(R^2 - 1)} \left( \rho - \dfrac{1}{\rho} \right) - \dfrac{\rho}{2} \log(\rho) \right]
\widehat{\theta} \\[6pt]
& q_{\lambda}(\xi) = \lambda^{2} \int_{1}^{\rho} \dfrac{1}{t} \left[ \dfrac{R^{2} \log(R)}{2(R^2 - 1)}  \left( t - \dfrac{1}{t} \right) -
\dfrac{t}{2} \log(t) \right]^{2} dt - \lambda \theta
\end{aligned}
\right. \qquad \forall \xi \in \overline{\Omega_{R}}
$$
is a nontrivial incomplete solution (with associated pressure) of the unforced Navier-Stokes system \eqref{nsstokes0} when $\omega=0$.
\end{remark}

\begin{remark}\label{alternative}
Consider a simple smooth line $\Sigma_{R}\subset\Omega_R$ such that $\Omega_R\setminus \overline{\Sigma_{R}}$ is simply connected. Following
\cite[Lemme 1.2]{foiastemam78}, let $\overline{q} \in \mathcal{C}^{\infty}(\Omega_R^*) \cap L^{\infty}(\Omega_R)$ be the unique scalar function
satisfying \eqref{firstchar}. Therefore, since $\overline{q}$ is discontinuous along $\Sigma_{R}$ (but its gradient can be extended continuously
to the whole $\Omega_{R}$), its distributional derivative involves a Dirac delta as in \eqref{dirac}. In this case, incomplete solutions to the
unforced Navier-Stokes system \eqref{nsstokes0} can be built in a similar way as in Theorem \ref{nonuniquetheo}, where the function
$\overline{q}$ induces a singularity on $\Sigma_R$ in the scalar pressure.
\end{remark}

The multiplicity result of Theorem \ref{nonuniquetheomain} can be partly explained because the difference of two incomplete solutions to
\eqref{nsstokes0} is not necessarily an element of $\mathbb{H}_{\sigma}(\Omega_{R})$. Nevertheless, the difference between two incomplete
solutions to \eqref{nsstokes0} having a \textit{prescribed} flux across $\Sigma_{R}$ indeed becomes a test function in $\mathbb{H}_{
\sigma}(\Omega_{R})$ and, therefore, one would expect unique solvability of \eqref{nsstokes0} in this class of incomplete solutions (namely,
with prescribed flux across $\Sigma_{R}$) under an appropriate smallness assumption on the data. With the notation as in \eqref{fluxcarrier},
given any $\omega \geq 0$ and $\Phi \in \mathbb{R}$, let us put
$$
\begin{aligned}
& K_{R}(\omega, \Phi) := \|\nabla V_{*}\|_{L^{2}(\Omega_{R})} = \sqrt{2 \pi} \Bigg[\dfrac{R^2 + 1}{R^2 - 1} \left( \omega + \dfrac{R \lambda_{*}}{2} \log(R) \right)^2 + \dfrac{\lambda^{2}_{*}}{8} (2 R^2 \log(R)^2 + R^2 - 1) \\[6pt]
& \hspace{5.3cm} - \dfrac{\lambda_{*} R(R^2 + 1) \log(R) }{R^2 - 1}\left( \omega + \dfrac{R \lambda_{*}}{2} \log(R) \right) \Bigg]^{1/2} \, .
\end{aligned}
$$

Defining $\mathcal{S}_4 >0$ as in \eqref{sobolevconstants1}, we now prove:

\begin{theorem} \label{fluxtheo1}
For any $\omega \geq 0$, $\Phi \in \mathbb{R}$ and $f \in \mathcal{R}[L^{2}(\Omega_{R})]$, there exist at least one rotationally-invariant
incomplete solution $V_{\Phi} \in \mathcal{R}[H^{1}_{*,\sigma}(\Omega_{R})]$ to \eqref{nsstokes0} having flux $\Phi$,
with associated pressure $Q_{\Phi} \in L^{2}(\Omega_{R})/\mathbb{R}$. Then, $(V_{\Phi},Q_{\Phi})$ is the unique incomplete solution to \eqref{nsstokes0} having flux $\Phi$
\begin{itemize}
\item in the space $\mathcal{R}[H^{1}_{*,\sigma}(\Omega_{R})] \times L^{2}(\Omega_{R})/\mathbb{R}$ for any $\omega \geq 0$, $f\in \mathcal{R}[L^{2}(\Omega_{R})]$, $\Phi \in \mathbb{R}$;
\item in the whole space $H^{1}_{*,\sigma}(\Omega_{R}) \times L^{2}(\Omega_{R})/\mathbb{R}$ whenever
\begin{equation} \label{umbral1cor3}
\dfrac{\|f\|_{L^2(\Omega_{R})}}{\sqrt{\lambda(\Omega_{R})}} + K_{R}(\omega, \Phi) < \mathcal{S}_4  \, .
\end{equation}
\end{itemize}
\end{theorem}
\noindent
\begin{proof} From the proof of Theorem
\ref{nonuniquetheomain} we know that there exists at least one rotationally-invariant incomplete solution $V_{\Phi} \in
\mathcal{R}[H^{1}_{*,\sigma}(\Omega_{R})]$ of \eqref{nsstokes0} having flux $\Phi$. Then, Theorem \ref{incompletepressure} ensures the existence
of a uniquely associated scalar pressure $Q_\Phi\in L^{2}(\Omega_{R})/\mathbb{R}$.\par
Concerning the first item of the statement, suppose that $V \in \mathcal{R}[H^{1}_{*,\sigma}(\Omega_{R})]$ is another incomplete solution to
\eqref{nsstokes0} having flux $\Phi$. Repeating the argument of the proof of Theorem \ref{nonuniquetheomain}, we deduce there exist two functions
$G_{\Phi}, G \in H^{1}(1,R)$ such that $G_{\Phi}(1)=G(1)=0$, $G_{\Phi}(R)=G(R)=\omega$ and
\begin{equation} \label{flux11}
V_{\Phi}(\xi) = G_{\Phi}(\rho) \widehat{\theta} \qquad \text{and} \qquad  V(\xi) = G(\rho) \widehat{\theta} \qquad \text{for a.e. }\xi \in
\Omega_{R} \, .
\end{equation}
Let $z := V - V_{\Phi} \in \mathbb{H}_{\sigma}(\Omega_{R})$ and subtract the equations \eqref{nstokesdebilstar} corresponding to $V$
and $V_{\Phi}$, thereby obtaining
$$
\int_{\Omega_{R}} \nabla z \cdot \nabla \varphi + \int_{\Omega_{R}} (V \cdot \nabla) z \cdot \varphi + \int_{\Omega_{R}} (z \cdot \nabla) V_{\Phi}
\cdot \varphi = 0 \qquad \forall \varphi \in \mathbb{H}_{\sigma}(\Omega_{R}) \, .
$$
By taking $\varphi=z$, integrating by parts, and using \eqref{flux11}, we get
$$
\int_{\Omega_{R}} (V \cdot \nabla) z \cdot z = 0 \qquad \text{and} \qquad \int_{\Omega_{R}} (z \cdot \nabla) V_{\Phi} \cdot z = -
\int_{\Omega_{R}} (z \cdot \nabla) z \cdot V_{\Phi} = 0 \, ,
$$
and we deduce that $z=0$ in $\Omega_{R}$. Then, Theorem \ref{incompletepressure} ensures that the associated pressures coincide (up to an additive constant)
a.e.\ in $\Omega_{R}$ as well.
\par
Concerning the second item of the statement, by taking $\varphi = V_{\Phi} - V_{*} \in \mathbb{H}_{\sigma}(\Omega_{R})$ in both weak formulations
\eqref{nstokesdebilstar} and \eqref{lambdastarwf}, and noticing that
$$
(V_{\Phi} \cdot \nabla) V_{\Phi} \cdot (V_{\Phi} - V_{*}) =0 \quad \text{a.e. in} \ \Omega_{R} \, ,
$$
we deduce
$$
\| \nabla V_{\Phi} - \nabla V_{*} \|^{2}_{L^{2}(\Omega_{R})} = \int_{\Omega_{R}} f \cdot (V_{\Phi} - V_{*}) \leq
\dfrac{\|f\|_{L^2(\Omega_{R})}}{\sqrt{\lambda(\Omega_{R})}} \| \nabla V_{\Phi} - \nabla V_{*} \|_{L^{2}(\Omega_{R})} \, ,
$$
as a consequence of H\"older's inequality and \eqref{sobpoin}. Therefore,
\begin{equation} \label{flux2}
\| \nabla V_{\Phi} \|_{L^{2}(\Omega_{R})} \leq \| \nabla V_{\Phi} - \nabla V_{*} \|_{L^{2}(\Omega_{R})} + \| \nabla V_{*} \|_{L^{2}(\Omega_{R})}
\leq \dfrac{\|f\|_{L^2(\Omega_{R})}}{\sqrt{\lambda(\Omega_{R})}} + \| \nabla V_{*} \|_{L^{2}(\Omega_{R})} \, .
\end{equation}
Now, suppose that $V \in H^{1}_{*,\sigma}(\Omega_{R})$ is another incomplete solution to  \eqref{nsstokes0} (not necessarily rotationally
invariant) having flux $\Phi$. Let $z := V - V_{\Phi} \in \mathbb{H}_{\sigma}(\Omega_{R})$ and subtract the equations
\eqref{nstokesdebilstar} corresponding to $V$ and $V_{\Phi}$, thereby obtaining
$$
\int_{\Omega_{R}} \nabla z \cdot \nabla \varphi + \int_{\Omega_{R}} (V \cdot \nabla) z \cdot \varphi + \int_{\Omega_{R}} (z \cdot \nabla) V_{\Phi}
\cdot \varphi = 0 \qquad \forall \varphi \in \mathbb{H}_{\sigma}(\Omega_{R}) \, .
$$
By taking $\varphi=z$ and noticing (after an integration by parts) that
$$
\int_{\Omega_{R}} (V \cdot \nabla) z \cdot z = 0 \, ,
$$
we deduce, from the H\"older inequality, from \eqref{sobolevconstants11} and \eqref{flux2}, that
$$
\|\nabla z\|^2_{L^2(\Omega_{R})} = - \int_{\Omega_{R}} (z \cdot \nabla)  V_{\Phi} \cdot z \leq \left(
\dfrac{\|f\|_{L^2(\Omega_{R})}}{\sqrt{\lambda(\Omega_{R})}} + \| \nabla V_{*} \|_{L^{2}(\Omega_{R})} \right) \frac{\|\nabla
z\|^2_{L^2(\Omega_{R})}}{\mathcal{S}_4}  \, ,
$$
so $z = 0$ a.e.\ in $\Omega_{R}$ whenever \eqref{umbral1cor3} holds.
\end{proof}

Theorem \ref{fluxtheo1} shows that if we enlarge the class of solutions by dropping one space dimension in the space of test functions,
see Definition \ref{incompletedef}, we may recover the same uniqueness result as in Theorem \ref{UVconstant0}. Roughly speaking, the additional
degree of freedom is deleted by the flux constraint.

\newpage
\section{Invading domains and a generalized Stokes paradox}\label{stokesparadox}
Let $\Omega_{\infty} := \mathbb{R}^2 \setminus \overline{B_1}$. The classical \textit{Stokes paradox} \cite{stokes1851effect} states that
the problem
\begin{equation}\label{stokesaradox}
	\left\{
	\begin{aligned}
		& -\Delta v+\nabla q=0 \, ,\ \quad  \nabla\cdot v=0 \ \ \mbox{ in } \ \ \Omega_{\infty} \, , \\[3pt]
		& \lim\limits_{\rho\to \infty} v(\xi) = (1,0) \, , \qquad v = (0,0) \ \text{ on } \ \partial B_{1} \, ,
	\end{aligned}
	\right.
\end{equation}
has no solution, see also \cite[Chapter V]{galdi2011introduction} and \cite{korobkov2023stationary}. Even more, Chang \& Finn
\cite{chang1961solutions} proved that the only solution to the exterior problem
\begin{equation}\label{stokesext}
	\left\{
	\begin{aligned}
		& -\Delta v+\nabla q=0 \, ,\ \quad  \nabla\cdot v=0 \ \ \mbox{ in } \ \ \Omega_{\infty} \, , \\[3pt]
		& v = (0,0) \ \text{ on } \ \partial B_{1} \, ,
	\end{aligned}
	\right.
\end{equation}
such that $|v(\xi)| = o(\log(\rho))$ as $\rho \to \infty$, is the trivial solution. Notice that \eqref{stokesaradox} may be obtained by
letting $R\to\infty$ in the annulus $\Omega_R$ (see \eqref{annulus}) for the steady-state Stokes equations
\begin{equation}\label{stokes0}
	\left\{
	\begin{aligned}
		& -\Delta v+\nabla q=0 \, ,\ \quad  \nabla\cdot v=0 \ \ \mbox{ in } \ \ \Omega_{R}, \\[3pt]
		& v = (1,0) \ \text{ on } \ \partial B_R \, , \qquad v = (0,0) \ \text{ on } \ \partial B_{1} \, ,
	\end{aligned}
	\right.
\end{equation}
that have an explicit solution. Indeed, consider the constants (which depend on $R>1$)
\begin{equation} \label{constantes}
	\begin{aligned}
		& C_1 := \dfrac{-1}{2(1-R^2 + \log(R) + R^2 \log(R))} \, , \qquad C_2 := \dfrac{R^2}{2(1-R^2 + \log(R) + R^2 \log(R))} \, ,\\[6pt]
		& C_3 := \dfrac{1-R^2}{2(1-R^2 + \log(R) + R^2 \log(R))} \, ,\qquad C_4 := \dfrac{R^2 + 1}{1-R^2 + \log(R) + R^2 \log(R)}\, .
	\end{aligned}
\end{equation}
Then, a direct computation proves the following result:
\begin{proposition} \label{exactstokes}
The unique (classical) solution $(v_R,q_R)\in\mathcal{C}^{2}(\overline{\Omega_{R}}) \times \mathcal{C}^{1}(\overline{\Omega_{R}})/\mathbb{R}$  of \eqref{stokes0} reads
	$$
	\left\{
	\begin{aligned}
		& v_{R}(\xi) =  \left( C_{1} \rho^{2} + \dfrac{C_{2}}{\rho^{2}} + C_{3} + C_{4} \log(\rho) \right) \cos(\theta) \widehat{\rho} -
\left( 3 C_{1} \rho^{2} - \dfrac{C_{2}}{\rho^{2}} + C_{3} + C_{4} (1+ \log(\rho)) \right) \sin(\theta) \widehat{\theta} \, , \\[6pt]
		& q_R(\xi) =   \left( 8C_1 \rho - \dfrac{2 C_4}{\rho} \right) \cos(\theta) \, ,
	\end{aligned}
	\right.
	$$
	where $C_1, C_2, C_3, C_4 \in \mathbb{R}$ are defined in
\eqref{constantes}.
\end{proposition}

From \eqref{constantes} we infer that, as $R\to\infty$,
\begin{equation}\label{todos0}
C_1(R) \to 0 \, ,\qquad C_2(R)\to0 \,,\qquad C_3(R)\to0\,,\qquad C_4(R)\to0\,,
\end{equation}
so that Proposition \ref{exactstokes} shows that ($v_R$ and $p_R$ being extended trivially in $\R^2\setminus B_R$)
$$
v_R, \, p_R\to0\quad\mbox{in }L^\infty_{\rm loc}(\Omega_\infty)\quad\mbox{as }R\to\infty\, .
$$
This is why the Leray
invading domains technique \cite{leray1933etude} does not allow to find solutions to \eqref{stokesaradox},
giving an interpretation to the Stokes paradox. From a physical point of view, this paradox can be explained as follows:
\begin{center}
a fixed inflow is too weak to maintain the movement of a highly viscous fluid\\
when the domain becomes too large.
\end{center}
The failure of the Leray invading domains technique, together with the appearance of the Stokes paradox, suggest to modify \eqref{stokes0}
by enlarging also the (inflow/outflow) boundary data and to consider the problem
\begin{equation}\label{stokes1}
\left\{
\begin{aligned}
	& -\Delta v+\nabla q=0,\ \quad  \nabla\cdot v=0 \ \ \mbox{ in } \ \ \Omega_{R}, \\[3pt]
	& v = (\log (R),0) \ \text{ on } \ \partial B_R, \qquad v = (0,0) \ \text{ on } \ \partial B_{1} \, .
\end{aligned}
\right.
\end{equation}

By linearity, we obtain the following straightforward consequence of Proposition \ref{exactstokes}:

\begin{corollary}\label{outerstokes}
The unique (classical) solution $(v^R,q^R)\in\mathcal{C}^{2}(\overline{\Omega_{R}}) \times \mathcal{C}^{1}(\overline{\Omega_{R}})/\mathbb{R}$ of \eqref{stokes1} reads
$$
\left\{
\begin{aligned}
& v^R(\xi)=\log(R)\left[ \left(\!C_{1} \rho^{2}\! +\! \dfrac{C_{2}}{\rho^{2}}\! +\! C_{3}\! +\! C_{4} \log(\rho)\!\right) \cos(\theta)
\widehat{\rho} -  \left(\!3 C_{1} \rho^{2}\! -\!\dfrac{C_{2}}{\rho^{2}}\!+\! C_{3}\! +\! C_{4} (1\!+\!\log(\rho))\!\right)
\sin(\theta)\widehat{\theta}\right],\\[6pt]
	& q^R(\xi) =   \log(R)\left( 8C_1 \rho - \dfrac{2 C_4}{\rho} \right) \cos(\theta) \, ,
\end{aligned}
\right.
$$
where $C_1, C_2, C_3, C_4 \in \mathbb{R}$ are defined in \eqref{constantes}.
\end{corollary}

The fundamental difference with \eqref{todos0} is that, as $R\to\infty$,
$$
C_1(R)\log(R)\to0 \, , \quad C_2(R)\log(R)\to\frac12 \, , \quad C_3(R)\log(R)\to-\frac12 \, , \quad C_4(R)\log(R)\to1 \, .
$$
Hence, the invading domains technique applied to \eqref{stokes1} yields

\begin{proposition}\label{explicit} The couple of functions
$$
\left\{
\begin{aligned}
	& v(\xi) = \left( \dfrac{1}{2\rho^{2}} -\frac12 + \log(\rho) \right) \cos(\theta) \widehat{\rho} + \left(\dfrac{1}{2\rho^{2}} -\frac12
	-\log(\rho)\right) \sin(\theta) \widehat{\theta} \qquad \forall \xi \in \Omega_{\infty} \, , \\[6pt]
	& q(\xi) =  -\dfrac{2}{\rho} \cos(\theta) \qquad \forall \xi \in \Omega_{\infty} \, ,
\end{aligned}
\right.
$$
is a (nontrivial) solution to the exterior problem \eqref{stokesext}.
\end{proposition}

The solution in Proposition \ref{explicit} satisfies $|v(\xi)|\asymp\log(\rho)$ as $\rho\to\infty$, showing that the nonexistence result of
\cite{chang1961solutions} is sharp and complementing the physical interpretation of the Stokes paradox:
\begin{center}
to maintain a visible movement of a highly viscous fluid when the domain becomes larger,\\
an inflow of increasing magnitude is needed which, at the limit, is not physically attainable.
\end{center}
With the Taylor-Couette boundary conditions \eqref{nsstokes0}$_2$, we generalize the Stokes paradox and we find a somehow surprising statement.
First of all, notice that the results in \cite{chang1961solutions} imply that the problem
\begin{equation}\label{stoextrot}
\left\{
\begin{aligned}
	& -\Delta v + \nabla q=0 \, ,\ \quad  \nabla\cdot v=0 \ \ \mbox{ in } \ \ \Omega_{\infty} \, , \\[3pt]
	& \lim\limits_{\rho \to \infty} v(\xi) = \widehat{\theta}\, , \qquad v = (0,0) \ \text{ on } \ \partial B_{1} \,
\end{aligned}
\right.
\end{equation}
has no solution. Then, in the annulus $\Omega_R$ (see \eqref{annulus}), consider the Stokes equations
\begin{equation}\label{stokesV}
\left\{
\begin{aligned}
	& -\Delta v+\nabla q=0 \, , \ \quad  \nabla\cdot v=0 \ \ \mbox{ in } \ \ \Omega_{R} \, , \\[3pt]
	& v = \widehat{\theta} \ \text{ on } \ \partial B_R \, , \qquad v = (0,0) \ \text{ on } \ \partial B_{1} \, ,
\end{aligned}
\right.
\end{equation}
which is the linear version of \eqref{nsstokes0}. The counterpart of Proposition \ref{exactstokes} reads:

\begin{proposition} \label{exactstokestheta}
The unique (classical) solution $(V_R,Q_R)\in\mathcal{C}^{2}(\overline{\Omega_{R}}) \times \mathcal{C}^{1}(\overline{\Omega_{R}})/\mathbb{R}$ to \eqref{stokesV} reads
\begin{equation} \label{tcs1}
	V_{R}(\xi) := \dfrac{R}{R^2 - 1} \left( \rho - \dfrac{1}{\rho} \right) \widehat{\theta} \qquad \text{and} \qquad Q_{R}(\xi) := 0
\qquad \forall \xi \in \Omega_{R} \, .
\end{equation}
\end{proposition}

Again, as $R \to \infty$, the solution \eqref{tcs1}, when trivially extended by $\widehat{\theta}$ to $\Omega_{\infty}$, locally converges
uniformly to zero, in line with what was observed for \eqref{stokes0}: the invading domains method produces a zero solution to
\eqref{stoextrot}.
In order to find a nontrivial solution as $R\to\infty$, instead of \eqref{stoextrot} we consider the problem
\begin{equation}\label{stokes2}
\left\{
\begin{aligned}
	& -\Delta v+\nabla q=0,\ \quad  \nabla\cdot v=0 \ \ \mbox{ in } \ \ \Omega_{R}, \\[3pt]
	& v = R\, \widehat{\theta} \ \text{ on } \ \partial B_R, \qquad v = (0,0) \ \text{ on } \ \partial B_{1} \, .
\end{aligned}
\right.
\end{equation}

By linearity, we obtain the following straightforward consequence of Proposition \ref{exactstokestheta}:

\begin{corollary}\label{outerstokes2}
The unique (classical) solution $(V^R,Q^R)\in\mathcal{C}^{2}(\overline{\Omega_{R}}) \times \mathcal{C}^{1}(\overline{\Omega_{R}})/\mathbb{R}$ to \eqref{stokes2} reads
$$
V^R(\xi) := \dfrac{R^2}{R^2-1}\left(\rho-\dfrac{1}{\rho}\right)\widehat{\theta}\qquad\text{and}\qquad
Q^R(\xi) := 0\qquad\forall\xi\in\Omega_{R} \, .
$$
\end{corollary}

Hence, in the case of Taylor-Couette boundary conditions, the invading domains technique yields
$$
V^\infty(\xi) := \left( \rho - \dfrac{1}{\rho} \right) \widehat{\theta} \qquad \text{and} \qquad Q^\infty(\xi) := 0 \qquad \forall \xi \in
\Omega_\infty\, ,
$$
as a (nontrivial) solution to the exterior problem \eqref{stokesext} which satisfies $|V_\infty(\xi)|\asymp\rho$
as $\rho\to\infty$. This solution should be compared with \eqref{uniqueK}. Finally, going back to the flux carrier \eqref{fluxcarrier}, a direct computation shows that the pair
$$
U^\infty(\xi) := \dfrac{\rho}{2} \log(\rho) \widehat{\theta} \qquad \text{and} \qquad P^\infty(\xi) := \theta \qquad \forall \xi \in
\Omega_\infty\, ,
$$
is a solution to the exterior problem \eqref{stokesext} in $\Omega_{\infty} \setminus \overline{\Sigma_{\infty}}$ verifying $|V_\infty(\xi)|\asymp \rho \log(\rho)$ as $\rho\to\infty$,
where
$$
\Sigma_{\infty} := \{ \xi \in \Omega_{\infty} \ | \ \theta = 0 \} = \{ (x,y) \in \Omega_{\infty} \ | \ x > 0 \, , \ y=0\} \, .
$$

The results in the present section connect the original Stokes paradox with the Taylor-Couette problem and suggest the following open questions:
\begin{itemize}
\item[(I)] prove/disprove that the only solutions to \eqref{stokesext} satisfying $|v(\xi)|\asymp\log(\rho)$ as $\rho\to\infty$ are such that $v(\xi)/\log(\rho) \to U\in\R^2\setminus\{0\}$ as $\rho\to\infty$;
\item[(II)] prove/disprove that there are no incomplete solutions (that are not smooth solutions) to \eqref{stokesext} satisfying $|v(\xi)|=o(\rho\log(\rho))$ as $\rho\to\infty$.
\end{itemize}

\appendix
\section{Appendix: bounds for Poincaré-Sobolev constants in the annulus}\label{bondsSob}

We denote the Sobolev constant of the embedding $H_0^{1}(\Omega_{R}) \subset L^{p}(\Omega_{R})$, for every $p \in [2,+\infty)$, by
\begin{equation} \label{sobolevconstants1}
\mathcal{S}_{p} := \min_{v \in H_0^{1}(\Omega_{R}) \setminus \{0\}} \ \frac{\|\nabla v\|^2_{L^2(\Omega_{R})}}{\|v\|^{2}_{L^{p}(\Omega_{R})}}
\qquad \forall p \in [2,+\infty) \, ,
\end{equation}
so that $\mathcal{S}_{p} > 0$ and
\begin{equation} \label{sobolevconstants11}
\mathcal{S}_{p} \, \|v\|^{2}_{L^{p}(\Omega_{R})} \leq \|\nabla v\|^2_{L^2(\Omega_{R})} \qquad \forall v \in H_0^{1}(\Omega_{R}) \, .
\end{equation}
Inequality \eqref{sobolevconstants11} is also valid for vector functions (with the same constant): if
$v = (v_1, v_2) \in H_{0}^{1}(\Omega_{R})$ is a vector field, by the Minkowski inequality we get
\begin{equation} \label{vectorornot}
\begin{aligned}
\|v\|^{p}_{L^p(\Omega_{R})} & = \left\| \, | v_{1} |^{2} + | v_{2} |^{2} \, \right\|^{p/2}_{L^{p/2}(\Omega_{R})} \leq \left(
\|v_{1}\|^{2}_{L^p(\Omega_{R})} + \|v_{2}\|^{2}_{L^p(\Omega_{R})} \right)^{p/2} \\[3pt]
& \leq \left( \dfrac{1}{\mathcal{S}_{p}} \right)^{p/2} \left( \| \nabla v_{1}\|^{2}_{L^2(\Omega_{R})} + \| \nabla v_{2}\|^{2}_{L^2(\Omega_{R})}
\right)^{p/2} =
\left( \dfrac{1}{\mathcal{S}_{p}} \right)^{p/2} \| \nabla v \|^{p}_{L^2(\Omega_{R})} \, .
\end{aligned}
\end{equation}
Notice that, in particular, $\mathcal{S}_{p}$ is a function of $R>1$. We set
\begin{equation} \label{eigen1}
\lambda(\Omega_{R}) := \mathcal{S}_{2}(R) \qquad \forall R>1 \, ,
\end{equation}
corresponding to the first eigenvalue of the Laplace-Dirichlet operator in $\Omega_{R}$, so that the Poincaré inequality in $\Omega_{R}$ (when
$p=2$ in \eqref{sobolevconstants11}) reads
\begin{equation} \label{sobpoin}
\|v\|_{L^{2}(\Omega_{R})} \leq \dfrac{1}{\sqrt{\lambda(\Omega_{R})}} \|\nabla v\|_{L^2(\Omega_{R})} \qquad \forall v \in H_0^{1}(\Omega_{R}) \, .
\end{equation}

\begin{figure}[H]
	\begin{center}
		\includegraphics[scale=0.4]{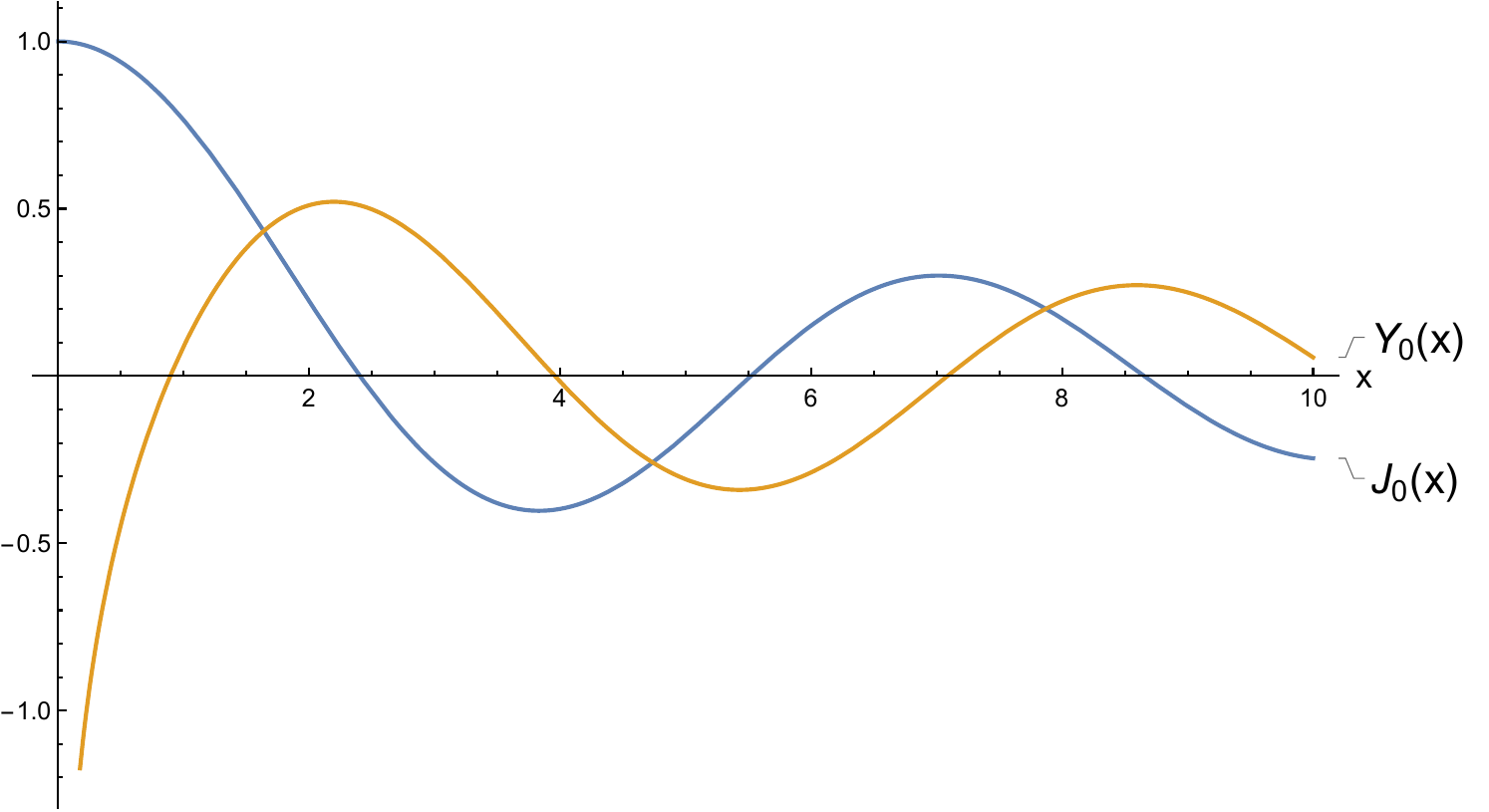}
	\end{center}
	\vspace*{-6mm}
	\caption{Plots of the Bessel functions $J_0$ and $Y_0$.}\label{bessel1}
\end{figure}

Let $J_{0} : [0,\infty) \longrightarrow \mathbb{R}$ and $Y_{0} : (0,\infty) \longrightarrow \mathbb{R}$ be, respectively, the Bessel functions of
the first and second kind of order zero, see \cite{AS} and
Figure \ref{bessel1}.

\begin{proposition}
Let $R>1$ and let $\alpha=\alpha(R)>0$ denote the first root of the equation
\begin{equation}\label{alfa}
Y_{0}(\alpha) J_{0} \left( \dfrac{\alpha}{R} \right) = J_{0}(\alpha) Y_{0} \left( \dfrac{\alpha}{R} \right) \, .
\end{equation}
Then,
\begin{equation}\label{lambda}
\dfrac{\pi^2}{2R^2}\le\lambda(\Omega_{R}) = \left( \dfrac{\alpha}{R} \right)^{2}\le \dfrac{10}{(R-1)^2}\, .
\end{equation}
\end{proposition}
\noindent
\begin{proof} The function
$$
F(\rho) := Y_{0}(\alpha) J_{0} \left( \dfrac{\alpha}{R} \rho \right) - J_{0}(\alpha) Y_{0} \left( \dfrac{\alpha}{R} \rho \right) \qquad
\forall \rho \in [1,R]
$$
vanishes both at $\rho =1$ and $\rho = R$ and has constant sign in $(1,R)$. Moreover, a simple computation shows that
$$
- \Delta F(\rho) = \left( \dfrac{\alpha}{R} \right)^{2} F(\rho) \qquad \forall \rho \in (1,R) \, ,
$$
implying the equality in \eqref{lambda}. The left inequality in \eqref{lambda} follows from the inclusion $\Omega_R\subset(-R,R)^2$ and the (decreasing) monotonicity of the
least Dirichlet eigenvalue of $-\Delta$: indeed, recall that in the square $(-R,R)^2$ it equals $\pi^2/2R^2$.
The right inequality in \eqref{lambda} follows by taking $v(\xi) := (\rho-1)(R-\rho)$ for $\xi \in \Omega_{R}$ and by computing the ratio in
\eqref{sobolevconstants1} for $p=2$.\end{proof}

\begin{center}
	\begin{minipage}{40mm}
		{\small
			\begin{tabular}{|c|c|c|}
				\hline
				Value of $R$ & Value of $\alpha$ & Value of $\lambda(\Omega_{R})$ \\
				\hline
				$ 1.1$ & $34.5535$ & $986.7308$ \\
				\hline
				$ 1.3$ & $13.6017$ & $109.4711$ \\
				\hline
				$ 1.5$ & $9.4053$ & $39.3158$ \\
				\hline
				$ 1.7$ & $7.6028$ & $20.0011$ \\
				\hline
				$2$ & $6.2461$ & $9.7533$ \\
				\hline
				$3$ & $4.6453$ & $2.3977$ \\
				\hline
				$4$ & $4.0976$ & $1.0494$ \\
				\hline
				$5$ & $3.8159$ & $0.5824$ \\
				\hline
				$6$ & $3.642$ & $0.3684$ \\
				\hline
				$7$ & $3.5227$ & $0.2532$ \\
				\hline
				$8$ & $3.4351$ & $0.1843$ \\
				\hline
				$9$ & $3.3677$ & $0.14002$ \\
				\hline
				$10$ & $3.3139$ & $0.1098$ \\
				\hline
				$15$ & $3.1504$ & $0.0441$ \\
				\hline
				$20$ & $3.0644$ & $0.0234$ \\
				\hline
				$100$ & $2.8009$ & $0.0007$ \\
				\hline
			\end{tabular}
		}
	\end{minipage}
	\qquad \qquad \qquad \qquad \qquad
	\begin{minipage}{85mm}
		\includegraphics[height=56mm,width=85mm]{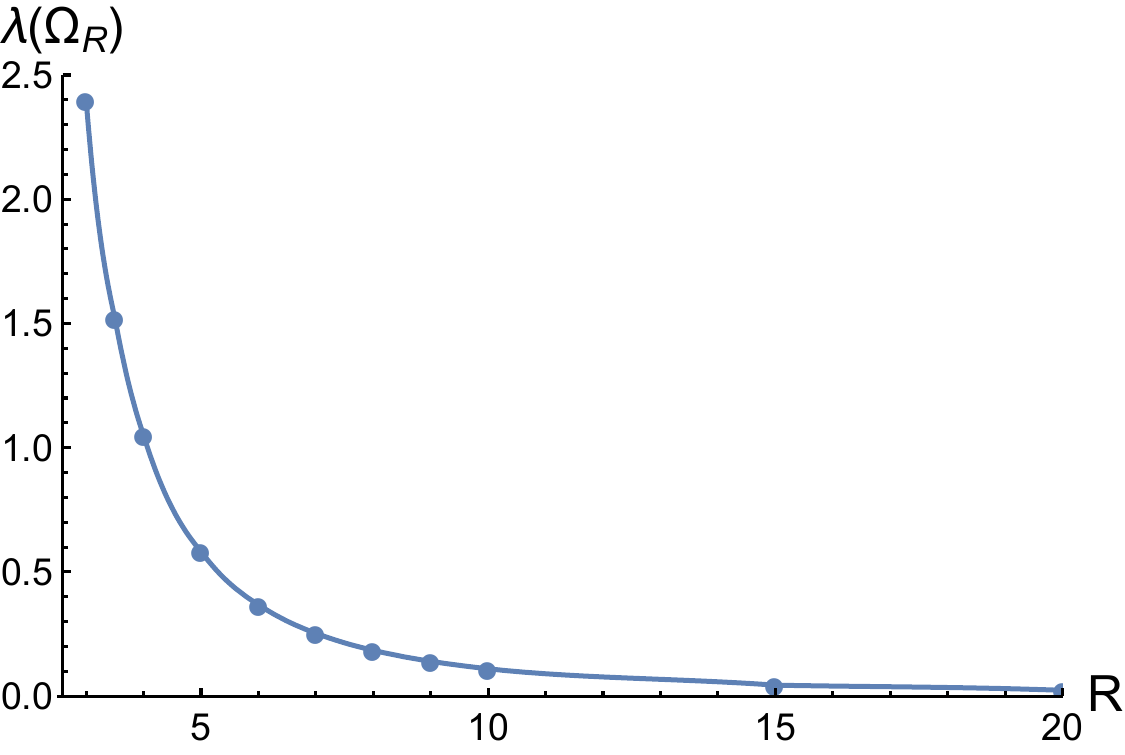}
	\end{minipage}
\end{center}
\begin{center}
	\vspace*{-3mm}	
	\captionof{figure}{Estimating the first Laplace-Dirichlet eigenvalue $\lambda(\Omega_{R})$ of $\Omega_{R}$ for different values of $R > 1$.}
	\label{figeigen}
\end{center}

It appears out of reach to solve equation \eqref{alfa} in closed form. Therefore, for different values of $R >1$, we give a
numerical approximation for the solution $\alpha$ of \eqref{alfa} and, in turn, for the eigenvalue $\lambda(\Omega_{R})$ in \eqref{lambda}.
This is complemented with a plot of the interpolating curve in Figure \ref{figeigen}.

We now turn our attention to bounds for a different Sobolev constant.
We introduce
\neweq{besselmu}
\mu_0=\mbox{the first zero of $J_0$}\approx2.40483\, ,
\endeq
and we prove

\begin{theorem}\label{boundsobolev}
Let $\Omega_R$ be as in \eqref{annulus}. For any scalar or vector function $w\in H^1_0(\Omega_{R})$ one has
	\neweq{ineqL4H1u2}
	\|w\|_{L^4(\Omega_{R})}^2 \leq \dfrac{1}{\mu_{0}} \sqrt{\dfrac{2}{3 \pi}(R^2 - 1)} \, \|\nabla w\|_{L^2(\Omega_{R})}^2\, .
	\endeq
In particular, for $\mathcal{S}_4$ as in \eqref{sobolevconstants1}, we have
\begin{equation} \label{lowerbound1}
\mu_{0} \sqrt{\dfrac{3 \pi}{2(R^2 - 1)}} \leq \mathcal{S}_4 \leq 24\sqrt{5 \pi} \, \dfrac{[(R^2 + 1)\log(R) + 1 - R^2](R^2 - 1)
\log(R)}{\sqrt{\kappa_{1}}} \, ,
\end{equation}
where we have defined
$$
\begin{aligned}
\kappa_{1} = & -1080 + 5400 R^2 - 10800 R^4 + 10800 R^6 - 5400 R^8 + 1080 R^{10} -
2025 \log(R) + 6075 R^2 \log(R) \\[3pt]
&  - 4050 R^4 \log(R) - 4050 R^6 \log(R) +
6075 R^8 \log(R) - 2025 R^{10} \log(R) - 1700 \log(R)^2  \\[3pt]
& + 3400 R^2 \log(R)^2 - 1700 R^4 \log(R)^2 + 1700 R^6 \log(R)^2 -
3400 R^8 \log(R)^2 + 1700 R^{10} \log(R)^2 \\[3pt]
& - 750 (1-R^2) \log(R)^3  + 750 R^8 \log(R)^3 - 750 R^{10} \log(R)^3 -
144 \log(R)^4 + 144 R^{10} \log(R)^4 .
\end{aligned}
$$
\end{theorem}
\noindent
\begin{proof} From \eqref{vectorornot} we know that it suffices to consider scalar functions $w \in H_{0}^{1}(\Omega_{R})$.
After combining the Gagliardo-Nirenberg inequality in $\R^2$ given by del Pino-Dolbeault \cite[Theorem 1]{delpino} with some H\"older inequality,
the following interpolation inequality was obtained in \cite[Theorem 2.3]{gazspe}:
\neweq{holder3}
\|w\|_{L^4(\Omega_{R})}^2\le \sqrt{\frac{2}{3\pi}} \, \|\nabla w\|_{L^2(\Omega_{R})}\|w\|_{L^2(\Omega_{R})}\qquad\forall w\in H^1_0(\Omega_R)\, ,
\endeq
which improves previous bounds by Ladyzhenskaya \cite{Lady59} (see also \cite[Lemma 1, p.8]{ladyzhenskaya1969mathematical}) and
\cite[Equation (II.3.9)]{galdi2011introduction}.		
Let $\Omega^* \subset \mathbb{R}^2$ be a disk having the same measure as $\Omega_R$, so that its radius is given by $\sqrt{R^2 - 1}$. Since the
Poincaré constant (least eigenvalue) in the unit disk is given by $\mu_0^2$, see \eqref{besselmu}, the Poincaré
constant of $\Omega^*$ is given by $\mu_0^2/(R^2-1)$. In view of the Faber-Krahn inequality \cite{faber1923beweis,krahn1925rayleigh} this means that
$$
\min_{w\in H^1_0(\Omega_R)}\ \frac{\|\nabla w\|_{L^2(\Omega_R)}}{\|w\|_{L^2(\Omega_R)}}\, \ge\, \min_{w\in H^1_0(\Omega^*)}\ \frac{\|\nabla
w\|_{L^2(\Omega^*)}}{\|w\|_{L^2(\Omega^*)}}
=\dfrac{\mu_0}{\sqrt{R^2 - 1}}\, .
$$
Therefore,
\begin{equation} \label{poinrum}
\|w\|_{L^2(\Omega_R)} \leq \frac{\sqrt{R^2 - 1}}{\mu_0}\|\nabla w\|_{L^2(\Omega_R)}
\qquad\forall w\in H^1_0(\Omega_R)
\end{equation}
which, once inserted into \eq{holder3}, yields \eq{ineqL4H1u2} and, consequently, the lower bound in \eqref{lowerbound1}. The upper bound in
\eqref{lowerbound1} follows by considering the function (defined in polar coordinates)
\begin{equation} \label{sobanubf}
X_{0}(\rho) := \dfrac{1}{4} \left[ (R^2 - 1) \dfrac{\log(\rho)}{\log(R)} + 1 - \rho^2 \right] \qquad \forall \rho \in [1,R] \, .
\end{equation}
\begin{figure}[H]
	\begin{center}
		\includegraphics[height=45mm,width=80mm]{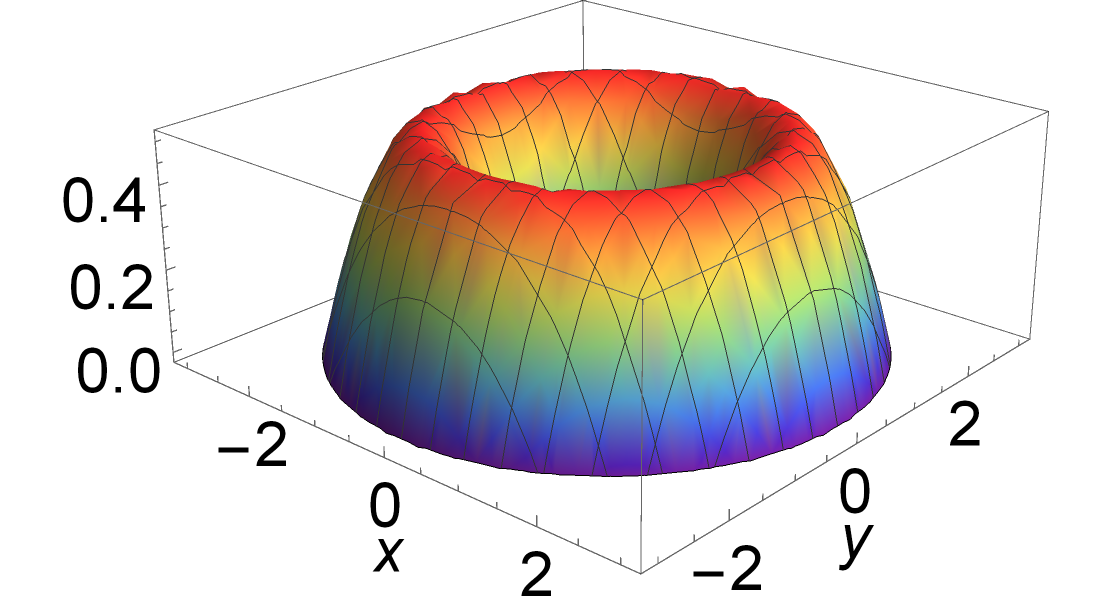} \qquad \includegraphics[height=45mm,width=80mm]{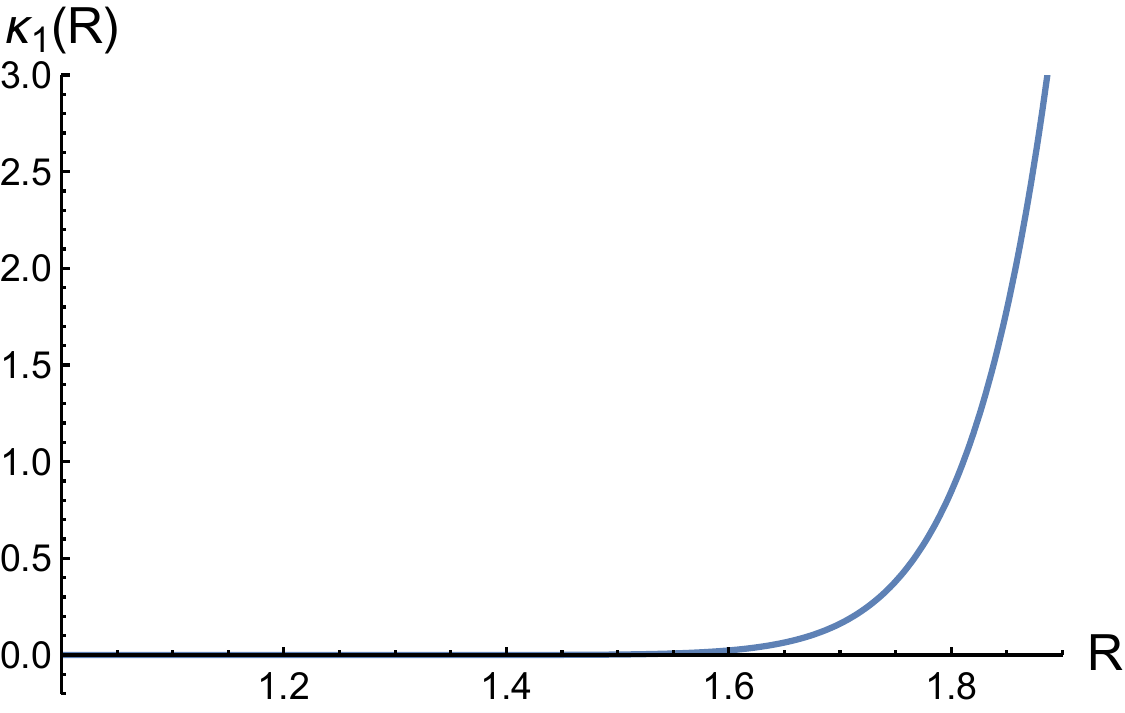}
	\end{center}
	\vspace*{-3mm}
	\caption{Left: Graph of $X_{0}$ in \eqref{sobanubf} for $R=3$. Right: plot of $\kappa_{1}$ in \eqref{lowerbound1} as a function of
$R>1$.}\label{sobanub}
\end{figure}
\noindent
Note that $X_0\in \mathcal{C}^\infty(\overline{\Omega_R})$ solves the torsion problem in $\Omega_R$:
$$
-\Delta X_{0} = 1 \ \ \mbox{ in } \ \ \Omega_R \, , \qquad X_{0}=0 \ \ \mbox{ on } \ \ \partial \Omega_R \, ,
$$
and can therefore be tested in the quotient \eqref{sobolevconstants1}.\end{proof}

\begin{remark} \label{remann}
From \eqref{lowerbound1} we observe that, as expected, $\mathcal{S}_4 \to 0$ as $R \to \infty$.
\end{remark}

We then turn to the task of estimating the optimal Sobolev embedding constant $\mathcal{S}_4$ defined in \eqref{sobolevconstants1}.
By combining \cite[Theorem 1.4]{lin1991} with \cite[Proposition 1.2]{nazarov2000one} (see also \cite[Theorem 5.3]{lin1992existence}),
we know that
\begin{equation}\label{yesno}
\mbox{the function achieving the minimum in \eqref{sobolevconstants1} is radial if $R\!\gg\!1$ and nonradial if $R\!\approx\!1$.}
\end{equation}
This suggests to introduce the subspace of $H_{0}^{1}(\Omega_{R})$ comprising radial functions:
$$
\mathcal{K}_{0}^{1}(\Omega_{R}) = \{v \in H_{0}^{1}(\Omega_{R}) \ | \ v(\xi) = v ( \rho ) \quad \forall \xi \in \Omega_R  \}\, ,
$$
and define the Sobolev constant of the embedding $\mathcal{K}_0^{1}(\Omega_{R}) \subset L^{4}(\Omega_{R})$ as
\begin{equation} \label{sobolevconstantrad}
\mathcal{R}_0 := \min_{v \in \mathcal{K}_{0}^{1}(\Omega_{R}) \setminus \{0\}} \ \dfrac{\|\nabla
v\|^2_{L^2(\Omega_{R})}}{\|v\|^{2}_{L^{4}(\Omega_{R})}} = \sqrt{2 \pi} \min_{v \in \mathcal{K}_{0}^{1}(\Omega_{R}) \setminus \{0\}} \
\dfrac{\displaystyle \int_{1}^{R} \rho |v'(\rho)|^{2} d\rho}{\displaystyle \sqrt{\int_{1}^{R} \rho |v(\rho)|^{4} d\rho}} \, .
\end{equation}

We can then rephrase \eq{yesno} as follows

\begin{corollary}\label{tobias}
For all $R > 1$ we have $\mathcal{S}_4 \leq \mathcal{R}_0$. Moreover, there exist $1<R_*\le R^*<\infty$ such that
\begin{center}
if \, $1<R<R_*$ \, then \, $\mathcal{S}_4<\mathcal{R}_0$\, ,\qquad \qquad if \, $R>R^*$ \, then \, $\mathcal{S}_4=\mathcal{R}_0$.
\end{center}
\end{corollary}

For any $R>1$, $\mathcal{R}_0$ provides an {\em upper bound} for $\mathcal{S}_4$ and this upper bound becomes an {\em equality}
if $R$ is sufficiently large, in which case the positive function achieving the minimum in \eqref{sobolevconstants1} is radial.\par

It is well-known \cite[Chapter I]{struwe2008variational} that any function $v \in H_{0}^{1}(\Omega_{R}) \setminus \{0\}$ achieving the minimum in
\eqref{sobolevconstants1} satisfies the following semilinear elliptic equation:
\begin{equation}\label{deq}
	\left\{
	\begin{aligned}
		& -\Delta v = v^3\ ,\quad  v > 0 \qquad \text{in} \quad \Omega_{R} \, , \\[5pt]
		& v = 0 \quad \text{on} \quad \partial\Omega_{R} \, .
	\end{aligned}
	\right.
\end{equation}
The advantage of restricting to radial functions $v\in \mathcal{K}_{0}^{1}(\Omega_{R})$ is that \eq{deq} becomes the ODE:
\begin{equation}\label{deqradial}
	\left\{
	\begin{aligned}
		& v''(\rho) + \dfrac{1}{\rho} v'(\rho) + v(\rho)^{3} = 0\, ,\quad v(\rho) > 0 \qquad \forall \rho \in (1,R) \, ,\\[5pt]
		& v(1) = v(R) = 0 \, .
	\end{aligned}
	\right.
\end{equation}

\begin{center}
	\begin{minipage}{5cm}
		\begin{tabular}{|c|c|c|}
			\hline
			Value of $R$ & Value of $a_{*}$ & Value of $\mathcal{R}_{0}$  \\
			\hline
			$ 1.1$ & $1004.0745$ & $644.592$  \\
			\hline
			$ 1.3$ & $118.258$ & $129.525$  \\
			\hline
			$ 1.5$ & $44.8402$ & $62.5316$  \\
			\hline
			$ 1.7$ & $23.9763$ & $39.057$  \\
			\hline
			$2$ & $12.5131$ & $23.9345$  \\
			\hline
			$3$ & $3.70134$ & $9.5259$  \\
			\hline
			$4$ & $1.87145$ & $5.6654$  \\
			\hline
			$5$ & $1.16988$ & $3.9517$  \\
			\hline
			$6$ & $0.8191$ & $3.00219$ \\
			\hline
			$7$ & $0.61524$ & $2.40486$  \\
			\hline
			$8$ & $0.48467$ & $1.99753$  \\
			\hline
			$9$ & $0.39518$ & $1.703$  \\
			\hline
			$10$ & $0.3307$ & $1.48117$   \\
			\hline
			$15$ & $0.17205$ & $0.8838$   \\
			\hline
			$20$ & $0.11095$ & $0.622396$  \\
			\hline
			$100$ & $0.012$ & $0.100268$  \\
			\hline
		\end{tabular}
	\end{minipage}
	\qquad \qquad \qquad \qquad
	\begin{minipage}{80mm}
		\includegraphics[height=54mm,width=81mm]{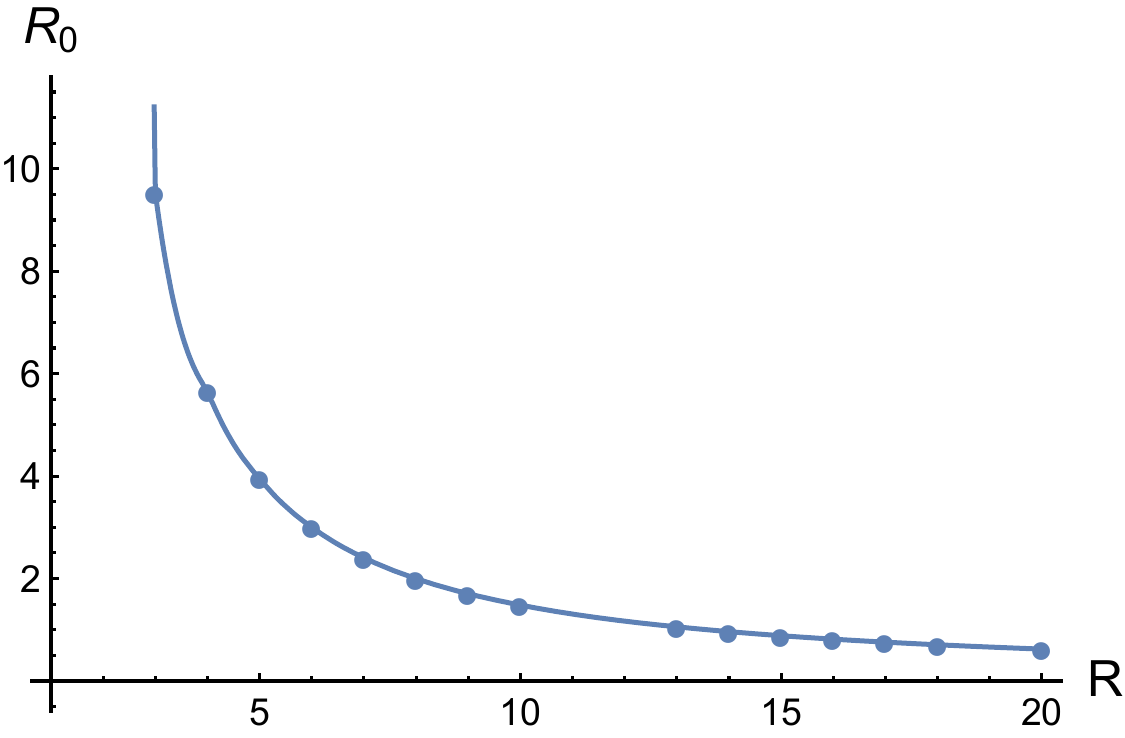}
	\end{minipage}
\end{center}
\begin{center}
	\vspace*{-3mm}	
	\captionof{figure}{Estimating the \textit{radial} Sobolev embedding constant \eqref{sobolevconstantrad} for different values of $R > 1$.}
	\label{radialsobb}	
\end{center}

We determine the solution to \eq{deqradial} numerically. In order to avoid the obtainment of the trivial solution $v \equiv 0$, we apply a \textit{shooting method}, that amounts to finding a solution to the problem
$$
\left\{
\begin{aligned}
& v''(\rho) + \dfrac{1}{\rho} v'(\rho) + v(\rho)^{3} = 0\, ,\quad v(\rho)>0\qquad \forall \rho \in (1,R) \, , \\[5pt]
& v(1) = 0 \, , \qquad v'(1) = a_{*} \, ,
\end{aligned}
\right.
$$
for varying values of $a_{*} > 0$, which is modified until a numerical solution to \eqref{deqradial} is found. With this procedure, for different
values of $R > 1$, we obtained an approximation of the embedding constant $\mathcal{R}_{0}$ given in \eqref{sobolevconstantrad}. Our results
are complemented with a plot of the interpolating curve in Figure \ref{radialsobb}. We are then interested in comparing the data from Figure \ref{radialsobb} with the analytic bounds given in \eqref{lowerbound1}.
This is displayed in Figure \ref{radialcomp} below.
\begin{figure}[H]
	\begin{center}
		\includegraphics[height=60mm,width=140mm]{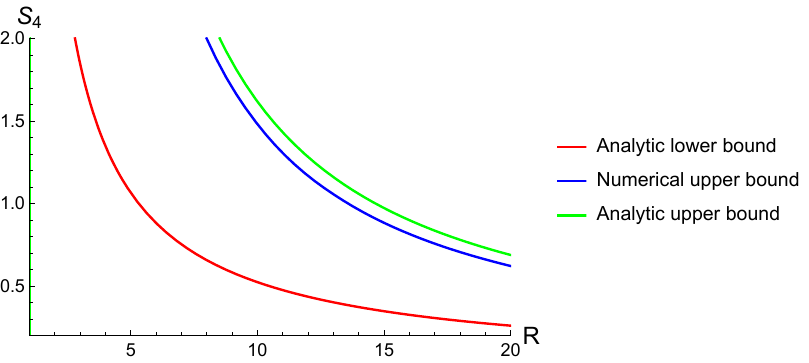}
	\end{center}
	\vspace*{-5mm}
	\caption{Behavior of the lower and upper bounds for $\mathcal{S}_4$ as functions of $R > 1$.}
	\label{radialcomp}
\end{figure}

As a mere curiosity, possibly suggesting conjectures on the behavior of $\lambda(\Omega_{R})$, let us notice:

\begin{remark}
By trivial extension and by the inclusions between $L^{q}(\Omega_{R})$ spaces we have, for any $q \in [2,+\infty)$, that
		\begin{equation}
		R\mapsto\mathcal{S}_q(R)\mbox{ is decreasing;} \qquad \qquad \mbox{if $R\le\sqrt{1+1/\pi}$, then }q\mapsto\mathcal{S}_q(R)\mbox{ is
decreasing,}
		\end{equation}
		while we expect $q\mapsto\mathcal{S}_q(R)$ to be increasing for large $R$. Moreover, by interpolation we have that
		\begin{eqnarray*}
			2\le p<q<r<\infty\ &\Longrightarrow&\ \|w\|^2_{L^q(\Omega_{R})}\le\|w\|^{2p(r-q)/q(r-p)}_{L^p(\Omega_{R})}
			\|w\|^{2r(q-p)/q(r-p)}_{L^r(\Omega_{R})}\qquad\forall w\in H^1_0(\Omega_R)\\[6pt]
			&\Longrightarrow&\ \mathcal{S}_q(R)\ge\mathcal{S}_p(R)^{p(r-q)/q(r-p)}\mathcal{S}_r(R)^{r(q-p)/q(r-p)} \, ,
		\end{eqnarray*}
		and also that
		$$
		2\le q\le\infty\ \Longrightarrow\ \|w\|^2_{L^q(\Omega_{R})}\le\|w\|^{4/q}_{L^2(\Omega_{R})}
		\|w\|^{2-4/q}_{L^\infty(\Omega_{R})}\qquad\forall w\in L^\infty(\Omega_R)\, .
		$$
\end{remark}

\bigskip\noindent
{\bf Acknowledgements.} The Department of Mathematics of the Politecnico di Milano was awarded by the MUR grant
{\em Dipartimento di Eccellenza 2023-27} (Italy) who supported this research. Most part of this paper was written while
Ji\v{r}\'{\i} Neustupa was visiting the Politecnico and Gianmarco Sperone was Assistant Professor at the Politecnico di Milano.
The work of Filippo Gazzola and Gianmarco Sperone is part of the PRIN project 2022 ``Partial differential equations and related
geometric-functional inequalities", financially supported by the EU, in the framework of the
``Next Generation EU initiative". Filippo Gazzola is also partially supported by INdAM. The work of Ji\v{r}\'{\i} Neustupa has been supported by the
Grant Agency of the Czech Republic (grant No.~22-01591S) and by the
Academy of Sciences of the Czech Republic (RVO 67985840). The Authors warmly thank the anonymous Referee for
careful proof-reading and for several suggestions that led to an improvement of the present work.\par
The Authors declare that there is no conflict of interest. Data sharing not applicable to this article as no datasets were generated or analyzed during the current study.

\newpage
\phantomsection
\addcontentsline{toc}{section}{References}
\bibliographystyle{abbrv}
\bibliography{references}
\vspace{5mm}
{\small
\begin{minipage}{56mm}
	Filippo Gazzola\\
	Dipartimento di Matematica\\
	Politecnico di Milano\\
	Piazza Leonardo da Vinci 32\\
	20133 Milan - Italy\\
	filippo.gazzola@polimi.it
\end{minipage}
\begin{minipage}{56mm}
	Ji\v{r}\'{\i} Neustupa\\
	Institute of Mathematics\\
	Czech Academy of Sciences\\
	\v{Z}itn\'{a} 25\\
	115 67 Prague - Czech Republic\\
	neustupa@math.cas.cz
\end{minipage}
\begin{minipage}{65mm}
	Gianmarco Sperone\\
Facultad de Matemáticas\\
Pontificia Universidad Católica de Chile\\
Avenida Vicuña Mackenna 4860\\
7820436 Santiago - Chile\\
E-mail: gianmarco.sperone@uc.cl
\end{minipage}
}
\end{document}